\documentclass{amsart}
\makeatletter
 
\def\diagram{\m@th\leftwidth=\z@ \rightwidth=\z@ \topheight=\z@
\botheight=\z@ \setbox\@picbox\hbox\bgroup}
 
\def\enddiagram{\egroup\wd\@picbox\rightwidth\unitlength
\ht\@picbox\topheight\unitlength \dp\@picbox\botheight\unitlength
\hskip\leftwidth\unitlength\box\@picbox}
 
\def\bfig{\begin{diagram}}
\def\efig{\end{diagram}}
\newcount\wideness \newcount\leftwidth \newcount\rightwidth
\newcount\highness \newcount\topheight \newcount\botheight
 
\def\ratchet#1#2{\ifnum#1<#2 \global #1=#2 \fi}
 
\def\putbox(#1,#2)#3{%
\horsize{\wideness}{#3} \divide\wideness by 2
{\advance\wideness by #1 \ratchet{\rightwidth}{\wideness}}
{\advance\wideness by -#1 \ratchet{\leftwidth}{\wideness}}
\vertsize{\highness}{#3} \divide\highness by 2
{\advance\highness by #2 \ratchet{\topheight}{\highness}}
{\advance\highness by -#2 \ratchet{\botheight}{\highness}}
\put(#1,#2){\makebox(0,0){$#3$}}}
 
\def\putlbox(#1,#2)#3{%
\horsize{\wideness}{#3}
{\advance\wideness by #1 \ratchet{\rightwidth}{\wideness}}
{\ratchet{\leftwidth}{-#1}}
\vertsize{\highness}{#3} \divide\highness by 2
{\advance\highness by #2 \ratchet{\topheight}{\highness}}
{\advance\highness by -#2 \ratchet{\botheight}{\highness}}
\put(#1,#2){\makebox(0,0)[l]{$#3$}}}
 
\def\putrbox(#1,#2)#3{%
\horsize{\wideness}{#3}
{\ratchet{\rightwidth}{#1}}
{\advance\wideness by -#1 \ratchet{\leftwidth}{\wideness}}
\vertsize{\highness}{#3} \divide\highness by 2
{\advance\highness by #2 \ratchet{\topheight}{\highness}}
{\advance\highness by -#2 \ratchet{\botheight}{\highness}}
\put(#1,#2){\makebox(0,0)[r]{$#3$}}}

\def\adjust[#1]{} % For compatibility
 
\newcount \coefa
\newcount \coefb
\newcount \coefc
\newcount\tempcounta
\newcount\tempcountb
\newcount\tempcountc
\newcount\tempcountd
\newcount\xext
\newcount\yext
\newcount\xoff
\newcount\yoff
\newcount\gap%
\newcount\arrowtypea
\newcount\arrowtypeb
\newcount\arrowtypec
\newcount\arrowtyped
\newcount\arrowtypee
\newcount\height
\newcount\width
\newcount\xpos
\newcount\ypos
\newcount\run
\newcount\rise
\newcount\arrowlength
\newcount\halflength
\newcount\arrowtype
\newdimen\tempdimen
\newdimen\xlen
\newdimen\ylen
\newsavebox{\tempboxa}%
\newsavebox{\tempboxb}%
\newsavebox{\tempboxc}%
 
\newdimen\w@dth
 
\def\setw@dth#1#2{\setbox\z@\hbox{\m@th$#1$}\w@dth=\wd\z@
\setbox\@ne\hbox{\m@th$#2$}\ifnum\w@dth<\wd\@ne \w@dth=\wd\@ne \fi
\advance\w@dth by 1.2em}
 
\def\t@^#1_#2{\allowbreak\def\n@one{#1}\def\n@two{#2}\mathrel
{\setw@dth{#1}{#2}
\mathop{\hbox to \w@dth{\rightarrowfill}}\limits
\ifx\n@one\empty\else ^{\box\z@}\fi
\ifx\n@two\empty\else _{\box\@ne}\fi}}
\def\t@@^#1{\@ifnextchar_{\t@^{#1}}{\t@^{#1}_{}}}
\def\to{\@ifnextchar^{\t@@}{\t@@^{}}}
 
\def\t@left^#1_#2{\def\n@one{#1}\def\n@two{#2}\mathrel{\setw@dth{#1}{#2}
\mathop{\hbox to \w@dth{\leftarrowfill}}\limits
\ifx\n@one\empty\else ^{\box\z@}\fi
\ifx\n@two\empty\else _{\box\@ne}\fi}}
\def\t@@left^#1{\@ifnextchar_{\t@left^{#1}}{\t@left^{#1}_{}}}
\def\toleft{\@ifnextchar^{\t@@left}{\t@@left^{}}}
 
\def\two@^#1_#2{\allowbreak
\def\n@one{#1}\def\n@two{#2}\mathrel{\setw@dth{#1}{#2}
\mathop{\vcenter{\lineskip\z@\baselineskip\z@
                 \hbox to \w@dth{\rightarrowfill}%
                 \hbox to \w@dth{\rightarrowfill}}%
       }\limits
\ifx\n@one\empty\else ^{\box\z@}\fi
\ifx\n@two\empty\else _{\box\@ne}\fi}}
\def\tw@@^#1{\@ifnextchar _{\two@^{#1}}{\two@^{#1}_{}}}
\def\two{\@ifnextchar ^{\tw@@}{\tw@@^{}}}
 
\def\tofr@^#1_#2{\def\n@one{#1}\def\n@two{#2}\mathrel{\setw@dth{#1}{#2}
\mathop{\vcenter{\hbox to \w@dth{\rightarrowfill}\kern-1.7ex
                 \hbox to \w@dth{\leftarrowfill}}%
       }\limits
\ifx\n@one\empty\else ^{\box\z@}\fi
\ifx\n@two\empty\else _{\box\@ne}\fi}}
\def\t@fr@^#1{\@ifnextchar_ {\tofr@^{#1}}{\tofr@^{#1}_{}}}
\def\tofro{\@ifnextchar^ {\t@fr@}{\t@fr@^{}}}

\def\mon{\mathop{\m@th\hbox to
      14.6\P@{\lasyb\char'51\hskip-2.1\P@$\arrext$\hss
$\mathord\rightarrow$}}\limits} % width of \epi
\def\leftmono{\mathrel{\m@th\hbox to
14.6\P@{$\mathord\leftarrow$\hss$\arrext$\hskip-2.1\P@\lasyb\char'50%
}}\limits} % width of \epi
\mathchardef\arrext="0200       % amr minus for arrow extension (see \into)

\def\settypes(#1,#2,#3){\arrowtypea#1 \arrowtypeb#2 \arrowtypec#3}
\def\settoheight#1#2{\setbox\@tempboxa\hbox{#2}#1\ht\@tempboxa\relax}%
\def\settodepth#1#2{\setbox\@tempboxa\hbox{#2}#1\dp\@tempboxa\relax}%
\def\settokens`#1`#2`#3`#4`{%
     \def\tokena{#1}\def\tokenb{#2}\def\tokenc{#3}\def\tokend{#4}}
\def\setsqparms[#1`#2`#3`#4;#5`#6]{%
\arrowtypea #1
\arrowtypeb #2
\arrowtypec #3
\arrowtyped #4
\width #5
\height #6
}
\def\setpos(#1,#2){\xpos=#1 \ypos#2}

\def\settriparms[#1`#2`#3;#4]{\settripairparms[#1`#2`#3`1`1;#4]}%
 
\def\settripairparms[#1`#2`#3`#4`#5;#6]{%
\arrowtypea #1
\arrowtypeb #2
\arrowtypec #3
\arrowtyped #4
\arrowtypee #5
\width #6
\height #6
}
 
\def\resetparms{\settripairparms[1`1`1`1`1;500]\width 500}%default values%
 
\resetparms
 
\def\mvector(#1,#2)#3{%%
\put(0,0){\vector(#1,#2){#3}}%
\put(0,0){\vector(#1,#2){26}}%
}
\def\evector(#1,#2)#3{{%%
\arrowlength #3
\put(0,0){\vector(#1,#2){\arrowlength}}%
\advance \arrowlength by-30
\put(0,0){\vector(#1,#2){\arrowlength}}%
}}
 
\def\horsize#1#2{%
\settowidth{\tempdimen}{$#2$}%
#1=\tempdimen
\divide #1 by\unitlength
}
 
\def\vertsize#1#2{%
\settoheight{\tempdimen}{$#2$}%
#1=\tempdimen
\settodepth{\tempdimen}{$#2$}%
\advance #1 by\tempdimen
\divide #1 by\unitlength
}
 
\def\putvector(#1,#2)(#3,#4)#5#6{{%
\ifnum3<\arrowtype
\putdashvector(#1,#2)(#3,#4)#5\arrowtype
\else
\ifnum\arrowtype<-3
\putdashvector(#1,#2)(#3,#4)#5\arrowtype
\else
\xpos=#1
\ypos=#2
\run=#3
\rise=#4
\arrowlength=#5
\ifnum \arrowtype<0
    \ifnum \run=0
        \advance \ypos by-\arrowlength
    \else
        \tempcounta \arrowlength
        \multiply \tempcounta by\rise
        \divide \tempcounta by\run
        \ifnum\run>0
            \advance \xpos by\arrowlength
            \advance \ypos by\tempcounta
        \else
            \advance \xpos by-\arrowlength
            \advance \ypos by-\tempcounta
        \fi
    \fi
    \multiply \arrowtype by-1
    \multiply \rise by-1
    \multiply \run by-1
\fi
\ifcase \arrowtype
\or \put(\xpos,\ypos){\vector(\run,\rise){\arrowlength}}%
\or \put(\xpos,\ypos){\mvector(\run,\rise)\arrowlength}%
\or \put(\xpos,\ypos){\evector(\run,\rise){\arrowlength}}%
\fi\fi\fi
}}
 
\def\putsplitvector(#1,#2)#3#4{%%
\xpos #1
\ypos #2
\arrowtype #4
\halflength #3
\arrowlength #3
\gap 140
\advance \halflength by-\gap
\divide \halflength by2
\ifnum\arrowtype>0
   \ifcase \arrowtype
   \or \put(\xpos,\ypos){\line(0,-1){\halflength}}%
       \advance\ypos by-\halflength
       \advance\ypos by-\gap
       \put(\xpos,\ypos){\vector(0,-1){\halflength}}%
   \or \put(\xpos,\ypos){\line(0,-1)\halflength}%
       \put(\xpos,\ypos){\vector(0,-1)3}%
       \advance\ypos by-\halflength
       \advance\ypos by-\gap
       \put(\xpos,\ypos){\vector(0,-1){\halflength}}%
   \or \put(\xpos,\ypos){\line(0,-1)\halflength}%
       \advance\ypos by-\halflength
       \advance\ypos by-\gap
       \put(\xpos,\ypos){\evector(0,-1){\halflength}}%
   \fi
\else \arrowtype=-\arrowtype
   \ifcase\arrowtype
   \or \advance \ypos by-\arrowlength
       \put(\xpos,\ypos){\line(0,1){\halflength}}%
       \advance\ypos by\halflength
       \advance\ypos by\gap
       \put(\xpos,\ypos){\vector(0,1){\halflength}}%
   \or \advance \ypos by-\arrowlength
       \put(\xpos,\ypos){\line(0,1)\halflength}%
       \put(\xpos,\ypos){\vector(0,1)3}%
       \advance\ypos by\halflength
       \advance\ypos by\gap
       \put(\xpos,\ypos){\vector(0,1){\halflength}}%
   \or \advance \ypos by-\arrowlength
       \put(\xpos,\ypos){\line(0,1)\halflength}%
       \advance\ypos by\halflength
       \advance\ypos by\gap
       \put(\xpos,\ypos){\evector(0,1){\halflength}}%
   \fi
\fi
}
 
\def\putmorphism(#1)(#2,#3)[#4`#5`#6]#7#8#9{{%
\run #2
\rise #3
\ifnum\rise=0
  \puthmorphism(#1)[#4`#5`#6]{#7}{#8}#9%
\else\ifnum\run=0
  \putvmorphism(#1)[#4`#5`#6]{#7}{#8}#9%
\else
\setpos(#1)%
\arrowlength #7
\arrowtype #8
\ifnum\run=0
\else\ifnum\rise=0
\else
\ifnum\run>0
    \coefa=1
\else
   \coefa=-1
\fi
\ifnum\arrowtype>0
   \coefb=0
   \coefc=-1
\else
   \coefb=\coefa
   \coefc=1
   \arrowtype=-\arrowtype
\fi
\width=2
\multiply \width by\run
\divide \width by\rise
\ifnum \width<0  \width=-\width\fi
\advance\width by60
\if l#9 \width=-\width\fi
\putbox(\xpos,\ypos){#4}%            %node 1
{\multiply \coefa by\arrowlength%      %node 2
\advance\xpos by\coefa
\multiply \coefa by\rise
\divide \coefa by\run
\advance \ypos by\coefa
\putbox(\xpos,\ypos){#5} }%
{\multiply \coefa by\arrowlength%      %label
\divide \coefa by2
\advance \xpos by\coefa
\advance \xpos by\width
\multiply \coefa by\rise
\divide \coefa by\run
\advance \ypos by\coefa
\if l#9%
   \putrbox(\xpos,\ypos){#6}%
\else\if r#9%
   \putlbox(\xpos,\ypos){#6}%
\fi\fi }%
{\multiply \rise by-\coefc%             %arrow
\multiply \run by-\coefc
\multiply \coefb by\arrowlength
\advance \xpos by\coefb
\multiply \coefb by\rise
\divide \coefb by\run
\advance \ypos by\coefb
\multiply \coefc by70
\advance \ypos by\coefc
\multiply \coefc by\run
\divide \coefc by\rise
\advance \xpos by\coefc
\multiply \coefa by140
\multiply \coefa by\run
\divide \coefa by\rise
\advance \arrowlength by\coefa
\ifcase\arrowtype
\or \put(\xpos,\ypos){\vector(\run,\rise){\arrowlength}}%
\or \put(\xpos,\ypos){\mvector(\run,\rise){\arrowlength}}%
\or \put(\xpos,\ypos){\evector(\run,\rise){\arrowlength}}%
\fi}\fi\fi\fi\fi}}

\newcount\numbdashes \newcount\lengthdash \newcount\increment
 
\def\howmanydashes{% Actually returns both number and length
\numbdashes=\arrowlength \lengthdash=40
\divide\numbdashes by \lengthdash
\lengthdash=\arrowlength
\divide\lengthdash by \numbdashes
\increment=\lengthdash
\multiply\lengthdash by 3
\divide\lengthdash by 5
}
 
\def\putdashvector(#1)(#2,#3)#4#5{%
\ifnum#3=0 \putdashhvector(#1){#4}#5
\else
\ifnum#2=0
\putdashvvector(#1){#4}#5\fi\fi}
 
\def\putdashhvector(#1,#2)#3#4{{%
\arrowlength=#3 \howmanydashes
\multiput(#1,#2)(\increment,0){\numbdashes}%
{\vrule height .4pt width \lengthdash\unitlength}
\arrowtype=#4 \xpos=#1
\ifnum\arrowtype<0 \advance\arrowtype by 7 \fi
\ifcase\arrowtype
\or \advance\xpos by 10
    \put(\xpos,#2){\vector(-1,0){\lengthdash}}
    \advance\xpos by 40
    \put(\xpos,#2){\vector(-1,0){\lengthdash}}
\or \advance \xpos by 10
    \put(\xpos,#2){\vector(-1,0){\lengthdash}}
    \advance\xpos by  \arrowlength
    \advance\xpos by  -50
    \put(\xpos,#2){\vector(-1,0){\lengthdash}}
\or \advance\xpos by 10
    \put(\xpos,#2){\vector(-1,0){\lengthdash}}
\or \advance\xpos by \arrowlength
    \advance\xpos by -\lengthdash
    \put(\xpos,#2){\vector(1,0){\lengthdash}}
\or {\advance\xpos by 10
    \put(\xpos,#2){\vector(1,0){\lengthdash}}}
    \advance\xpos by \arrowlength
    \advance\xpos by -\lengthdash
    \put(\xpos,#2){\vector(1,0){\lengthdash}}
\or \advance\xpos by \arrowlength
    \advance\xpos by -\lengthdash
    \put(\xpos,#2){\vector(1,0){\lengthdash}}
    \advance\xpos by -40
    \put(\xpos,#2){\vector(1,0){\lengthdash}}
   \fi
}}
 
\def\putdashvvector(#1,#2)#3#4{{%
\arrowlength=#3 \howmanydashes
\ypos=#2 \advance\ypos by -\arrowlength
\multiput(#1,#2)(0,\increment){\numbdashes}%
    {\vrule width .4pt height \lengthdash\unitlength}
\arrowtype=#4 \ypos=#2
\ifnum\arrowtype<0 \advance\arrowtype by 7 \fi
\ifcase\arrowtype
\or \advance\ypos by \arrowlength \advance\ypos by -40
    \put(#1,\ypos){\vector(0,1){\lengthdash}}
    \advance\ypos by -40
    \put(#1,\ypos){\vector(0,1){\lengthdash}}
\or \advance\ypos by 10
    \put(#1,\ypos){\vector(0,1){\lengthdash}}
    \advance\ypos by \arrowlength \advance\ypos by -40
    \put(#1,\ypos){\vector(0,1){\lengthdash}}
\or \advance\ypos by \arrowlength \advance\ypos by -40
    \put(#1,\ypos){\vector(0,1){\lengthdash}}
\or \advance\ypos by 10
    \put(#1,\ypos){\vector(0,-1){\lengthdash}}
\or \advance\ypos by 10
    \put(#1,\ypos){\vector(0,-1){\lengthdash}}
    \advance\ypos by \arrowlength \advance\ypos by -40
    \put(#1,\ypos){\vector(0,-1){\lengthdash}}
\or \advance\ypos by 10
    \put(#1,\ypos){\vector(0,-1){\lengthdash}}
    \advance\ypos by 40
    \put(#1,\ypos){\vector(0,-1){\lengthdash}}
\fi
}}
 
\def\puthmorphism(#1,#2)[#3`#4`#5]#6#7#8{{%
\xpos #1
\ypos #2
\width #6
\arrowlength #6
\arrowtype=#7
\putbox(\xpos,\ypos){#3\vphantom{#4}}%
{\advance \xpos by\arrowlength
\putbox(\xpos,\ypos){\vphantom{#3}#4}}%
\horsize{\tempcounta}{#3}%
\horsize{\tempcountb}{#4}%
\divide \tempcounta by2
\divide \tempcountb by2
\advance \tempcounta by30
\advance \tempcountb by30
\advance \xpos by\tempcounta
\advance \arrowlength by-\tempcounta
\advance \arrowlength by-\tempcountb
\putvector(\xpos,\ypos)(1,0)\arrowlength\arrowtype
\divide \arrowlength by2
\advance \xpos by\arrowlength
\vertsize{\tempcounta}{#5}%
\divide\tempcounta by2
\advance \tempcounta by20
\if a#8 %
   \advance \ypos by\tempcounta
   \putbox(\xpos,\ypos){#5}%
\else
   \advance \ypos by-\tempcounta
   \putbox(\xpos,\ypos){#5}%
\fi}}
 
\def\putvmorphism(#1,#2)[#3`#4`#5]#6#7#8{{%
\xpos #1
\ypos #2
\arrowlength #6
\arrowtype #7
\settowidth{\xlen}{$#5$}%
\putbox(\xpos,\ypos){#3}%
{\advance \ypos by-\arrowlength
\putbox(\xpos,\ypos){#4}}%
{\advance\arrowlength by-140
\advance \ypos by-70
\ifdim\xlen>0pt
   \if m#8%
      \putsplitvector(\xpos,\ypos)\arrowlength\arrowtype
   \else
   \putvector(\xpos,\ypos)(0,-1)\arrowlength\arrowtype
   \fi
\else
   \putvector(\xpos,\ypos)(0,-1)\arrowlength\arrowtype
\fi}%
\ifdim\xlen>0pt
   \divide \arrowlength by2
   \advance\ypos by-\arrowlength
   \if l#8%
      \advance \xpos by-40
      \putrbox(\xpos,\ypos){#5}%
   \else\if r#8%
      \advance \xpos by40
      \putlbox(\xpos,\ypos){#5}%
   \else
      \putbox(\xpos,\ypos){#5}%
   \fi\fi
\fi
}}
 
\def\putsquarep<#1>(#2)[#3;#4`#5`#6`#7]{{%
\setsqparms[#1]%
\setpos(#2)%
\settokens`#3`%
\puthmorphism(\xpos,\ypos)[\tokenc`\tokend`{#7}]{\width}{\arrowtyped}b%
\advance\ypos by \height
\puthmorphism(\xpos,\ypos)[\tokena`\tokenb`{#4}]{\width}{\arrowtypea}a%
\putvmorphism(\xpos,\ypos)[``{#5}]{\height}{\arrowtypeb}l%
\advance\xpos by \width
\putvmorphism(\xpos,\ypos)[``{#6}]{\height}{\arrowtypec}r%
}}
 
\def\putsquare{\@ifnextchar <{\putsquarep}{\putsquarep%
   <\arrowtypea`\arrowtypeb`\arrowtypec`\arrowtyped;\width`\height>}}
\def\square{\@ifnextchar< {\squarep}{\squarep
   <\arrowtypea`\arrowtypeb`\arrowtypec`\arrowtyped;\width`\height>}}
\def\squarep<#1>[#2`#3`#4`#5;#6`#7`#8`#9]{{%       %     #2------>#3
\setsqparms[#1]%                                   %      |       |
\diagram%                                          %      |       |
\putsquarep<\arrowtypea`\arrowtypeb`\arrowtypec`%  %    #7|       |#8
\arrowtyped;\width`\height>%                       %      |       |
(0,0)[#2`#3`#4`{#5};#6`#7`#8`{#9}]%                %      |       |
\enddiagram%                                       %      v       v
}}                                                 %     #4------>#5
\def\putptrianglep<#1>(#2,#3)[#4`#5`#6;#7`#8`#9]{{%
\settriparms[#1]%
\xpos=#2 \ypos=#3
\advance\ypos by \height
\puthmorphism(\xpos,\ypos)[#4`#5`{#7}]{\height}{\arrowtypea}a%
\putvmorphism(\xpos,\ypos)[`#6`{#8}]{\height}{\arrowtypeb}l%
\advance\xpos by\height
\putmorphism(\xpos,\ypos)(-1,-1)[``{#9}]{\height}{\arrowtypec}r%
}}
 
\def\putptriangle{\@ifnextchar <{\putptrianglep}{\putptrianglep
   <\arrowtypea`\arrowtypeb`\arrowtypec;\height>}}
\def\ptriangle{\@ifnextchar <{\ptrianglep}{\ptrianglep
   <\arrowtypea`\arrowtypeb`\arrowtypec;\height>}}
\def\ptrianglep<#1>[#2`#3`#4;#5`#6`#7]{{%%    %      #2----->#3
\settriparms[#1]%                             %      |      /
\diagram%                                     %      |     /
\putptrianglep<\arrowtypea`\arrowtypeb`%      %    #6|    /#7
\arrowtypec;\height>%                         %      |   /
(0,0)[#2`#3`#4;#5`#6`{#7}]%                   %      |  /
\enddiagram%%                                 %      v v
}}                                            %      #4
 
\def\putqtrianglep<#1>(#2,#3)[#4`#5`#6;#7`#8`#9]{{%
\settriparms[#1]%
\xpos=#2 \ypos=#3
\advance\ypos by\height
\puthmorphism(\xpos,\ypos)[#4`#5`{#7}]{\height}{\arrowtypea}a%
\putmorphism(\xpos,\ypos)(1,-1)[``{#8}]{\height}{\arrowtypeb}l%
\advance\xpos by\height
\putvmorphism(\xpos,\ypos)[`#6`{#9}]{\height}{\arrowtypec}r%
}}
 
\def\putqtriangle{\@ifnextchar <{\putqtrianglep}{\putqtrianglep
   <\arrowtypea`\arrowtypeb`\arrowtypec;\height>}}
\def\qtriangle{\@ifnextchar <{\qtrianglep}{\qtrianglep
   <\arrowtypea`\arrowtypeb`\arrowtypec;\height>}}
\def\qtrianglep<#1>[#2`#3`#4;#5`#6`#7]{{%%    %        #2----->#3
\settriparms[#1]%                             %         \      |
\width=\height                                %          \     |
\diagram%                                     %         #6\    |#7
\putqtrianglep<\arrowtypea`\arrowtypeb`%      %            \   |
\arrowtypec;\height>%                         %             \  |
(0,0)[#2`#3`#4;#5`#6`{#7}]%                   %              v v
\enddiagram%%                                 %               #4
}}
 
\def\putdtrianglep<#1>(#2,#3)[#4`#5`#6;#7`#8`#9]{{%
\settriparms[#1]%
\xpos=#2 \ypos=#3
\puthmorphism(\xpos,\ypos)[#5`#6`{#9}]{\height}{\arrowtypec}b%
\advance\xpos by \height \advance\ypos by\height
\putmorphism(\xpos,\ypos)(-1,-1)[``{#7}]{\height}{\arrowtypea}l%
\putvmorphism(\xpos,\ypos)[#4``{#8}]{\height}{\arrowtypeb}r%
}}
 
\def\putdtriangle{\@ifnextchar <{\putdtrianglep}{\putdtrianglep
   <\arrowtypea`\arrowtypeb`\arrowtypec;\height>}}
\def\dtriangle{\@ifnextchar <{\dtrianglep}{\dtrianglep
   <\arrowtypea`\arrowtypeb`\arrowtypec;\height>}}
\def\dtrianglep<#1>[#2`#3`#4;#5`#6`#7]{{%%    %                  / |
\settriparms[#1]%                             %                 /  |
\width=\height                                %              #5/   |#6
\diagram%                                     %               /    |
\putdtrianglep<\arrowtypea`\arrowtypeb`%      %              /     |
\arrowtypec;\height>%                         %             v      v
(0,0)[#2`#3`#4;#5`#6`{#7}]%                   %            #3----->#4
\enddiagram%%                                 %                #7
}}
 
\def\putbtrianglep<#1>(#2,#3)[#4`#5`#6;#7`#8`#9]{{%
\settriparms[#1]%
\xpos=#2 \ypos=#3
\puthmorphism(\xpos,\ypos)[#5`#6`{#9}]{\height}{\arrowtypec}b%
\advance\ypos by\height
\putmorphism(\xpos,\ypos)(1,-1)[``{#8}]{\height}{\arrowtypeb}r%
\putvmorphism(\xpos,\ypos)[#4``{#7}]{\height}{\arrowtypea}l%
}}
 
\def\putbtriangle{\@ifnextchar <{\putbtrianglep}{\putbtrianglep
   <\arrowtypea`\arrowtypeb`\arrowtypec;\height>}}
\def\btriangle{\@ifnextchar <{\btrianglep}{\btrianglep
   <\arrowtypea`\arrowtypeb`\arrowtypec;\height>}}
\def\btrianglep<#1>[#2`#3`#4;#5`#6`#7]{{%%   %              | \
\settriparms[#1]%                            %              |  \
\width=\height                               %            #5|   \#6
\diagram%                                    %              |    \
\putbtrianglep<\arrowtypea`\arrowtypeb`%     %              |     \
\arrowtypec;\height>%                        %              v      v
(0,0)[#2`#3`#4;#5`#6`{#7}]%                  %              #3----->#4
\enddiagram%%                                %                 #7
}}
 
\def\putAtrianglep<#1>(#2,#3)[#4`#5`#6;#7`#8`#9]{{%
\settriparms[#1]%
\xpos=#2 \ypos=#3
{\multiply \height by2
\puthmorphism(\xpos,\ypos)[#5`#6`{#9}]{\height}{\arrowtypec}b}%
\advance\xpos by\height \advance\ypos by\height
\putmorphism(\xpos,\ypos)(-1,-1)[#4``{#7}]{\height}{\arrowtypea}l%
\putmorphism(\xpos,\ypos)(1,-1)[``{#8}]{\height}{\arrowtypeb}r%
}}
 
\def\putAtriangle{\@ifnextchar <{\putAtrianglep}{\putAtrianglep
   <\arrowtypea`\arrowtypeb`\arrowtypec;\height>}}
\def\Atriangle{\@ifnextchar <{\Atrianglep}{\Atrianglep
   <\arrowtypea`\arrowtypeb`\arrowtypec;\height>}}
\def\Atrianglep<#1>[#2`#3`#4;#5`#6`#7]{{%%         %         /   \
\settriparms[#1]%                                  %        /     \
\width=\height                                     %     #5/       \#6
\diagram%                                          %      /         \
\putAtrianglep<\arrowtypea`\arrowtypeb`%           %     /           \
\arrowtypec;\height>%                              %    v             v
(0,0)[#2`#3`#4;#5`#6`{#7}]%                        %   #3------------>#4
\enddiagram%%                                      %          #7
}}
 
\def\putAtrianglepairp<#1>(#2)[#3;#4`#5`#6`#7`#8]{{%
\settripairparms[#1]%
\setpos(#2)%
\settokens`#3`%
\puthmorphism(\xpos,\ypos)[\tokenb`\tokenc`{#7}]{\height}{\arrowtyped}b%
\advance\xpos by\height
\puthmorphism(\xpos,\ypos)[\phantom{\tokenc}`\tokend`{#8}]%
{\height}{\arrowtypee}b%
\advance\ypos by\height
\putmorphism(\xpos,\ypos)(-1,-1)[\tokena``{#4}]{\height}{\arrowtypea}l%
\putvmorphism(\xpos,\ypos)[``{#5}]{\height}{\arrowtypeb}m%
\putmorphism(\xpos,\ypos)(1,-1)[``{#6}]{\height}{\arrowtypec}r%
}}
 
\def\putAtrianglepair{\@ifnextchar <{\putAtrianglepairp}{\putAtrianglepairp%
   <\arrowtypea`\arrowtypeb`\arrowtypec`\arrowtyped`\arrowtypee;\height>}}
\def\Atrianglepair{\@ifnextchar <{\Atrianglepairp}{\Atrianglepairp%
   <\arrowtypea`\arrowtypeb`\arrowtypec`\arrowtyped`\arrowtypee;\height>}}
 
\def\Atrianglepairp<#1>[#2;#3`#4`#5`#6`#7]{{%           %  #2a
\settripairparms[#1]%                         %           / | \
\settokens`#2`%                               %          /  |  \
\width=\height                                %       #3/  #4   \#5
\diagram%                                     %        /    |    \
\putAtrianglepairp                            %       /     |     \
<\arrowtypea`\arrowtypeb`\arrowtypec`%        %      v      v      v
\arrowtyped`\arrowtypee;\height>%             %     #2b---->#2c---->#2d
(0,0)[{#2};#3`#4`#5`#6`{#7}]%                 %         #6     #7
\enddiagram%%
}}
 
\def\putVtrianglep<#1>(#2,#3)[#4`#5`#6;#7`#8`#9]{{%
\settriparms[#1]%
\xpos=#2 \ypos=#3
\advance\ypos by\height
{\multiply\height by2
\puthmorphism(\xpos,\ypos)[#4`#5`{#7}]{\height}{\arrowtypea}a}%
\putmorphism(\xpos,\ypos)(1,-1)[`#6`{#8}]{\height}{\arrowtypeb}l%
\advance\xpos by\height
\advance\xpos by\height
\putmorphism(\xpos,\ypos)(-1,-1)[``{#9}]{\height}{\arrowtypec}r%
}}
 
\def\putVtriangle{\@ifnextchar <{\putVtrianglep}{\putVtrianglep
   <\arrowtypea`\arrowtypeb`\arrowtypec;\height>}}
\def\Vtriangle{\@ifnextchar <{\Vtrianglep}{\Vtrianglep
   <\arrowtypea`\arrowtypeb`\arrowtypec;\height>}}
\def\Vtrianglep<#1>[#2`#3`#4;#5`#6`#7]{{%%     %        #2------------->#3
\settriparms[#1]%                              %         \             /
\width=\height                                 %          \           /
\diagram%                                      %         #6\         /#7
\putVtrianglep<\arrowtypea`\arrowtypeb`%       %            \       /
\arrowtypec;\height>%                          %             \     /
(0,0)[#2`#3`#4;#5`#6`{#7}]%                    %              v   v
\enddiagram%%                                  %               #4
}}
 
\def\putVtrianglepairp<#1>(#2)[#3;#4`#5`#6`#7`#8]{{
\settripairparms[#1]%
\setpos(#2)%
\settokens`#3`%
\advance\ypos by\height
\putmorphism(\xpos,\ypos)(1,-1)[`\tokend`{#6}]{\height}{\arrowtypec}l%
\puthmorphism(\xpos,\ypos)[\tokena`\tokenb`{#4}]{\height}{\arrowtypea}a%
\advance\xpos by\height
\puthmorphism(\xpos,\ypos)[\phantom{\tokenb}`\tokenc`{#5}]%
{\height}{\arrowtypeb}a%
\putvmorphism(\xpos,\ypos)[``{#7}]{\height}{\arrowtyped}m%
\advance\xpos by\height
\putmorphism(\xpos,\ypos)(-1,-1)[``{#8}]{\height}{\arrowtypee}r%
}}
 
\def\putVtrianglepair{\@ifnextchar <{\putVtrianglepairp}{\putVtrianglepairp%
    <\arrowtypea`\arrowtypeb`\arrowtypec`\arrowtyped`\arrowtypee;\height>}}
\def\Vtrianglepair{\@ifnextchar <{\Vtrianglepairp}{\Vtrianglepairp%
    <\arrowtypea`\arrowtypeb`\arrowtypec`\arrowtyped`\arrowtypee;\height>}}
\def\Vtrianglepairp<#1>[#2;#3`#4`#5`#6`#7]{{%  %  #2a---->#2b---->#2c
\settripairparms[#1]%                          %   \      |      /
\settokens`#2`%                                %    \     |     /
\diagram%                                      %   #5\   #6    /#7
\putVtrianglepairp                             %      \   |   /
<\arrowtypea`\arrowtypeb`\arrowtypec`%         %       \  |  /
\arrowtyped`\arrowtypee;\height>%              %        v v v
(0,0)[{#2};#3`#4`#5`#6`{#7}]%                  %         #2d
\enddiagram%%
}}

\def\putCtrianglep<#1>(#2,#3)[#4`#5`#6;#7`#8`#9]{{%
\settriparms[#1]%
\xpos=#2 \ypos=#3
\advance\ypos by\height
\putmorphism(\xpos,\ypos)(1,-1)[``{#9}]{\height}{\arrowtypec}l%
\advance\xpos by\height
\advance\ypos by\height
\putmorphism(\xpos,\ypos)(-1,-1)[#4`#5`{#7}]{\height}{\arrowtypea}l%
{\multiply\height by 2
\putvmorphism(\xpos,\ypos)[`#6`{#8}]{\height}{\arrowtypeb}r}%
}}
 
\def\putCtriangle{\@ifnextchar <{\putCtrianglep}{\putCtrianglep
    <\arrowtypea`\arrowtypeb`\arrowtypec;\height>}}
\def\Ctriangle{\@ifnextchar <{\Ctrianglep}{\Ctrianglep
    <\arrowtypea`\arrowtypeb`\arrowtypec;\height>}}
\def\Ctrianglep<#1>[#2`#3`#4;#5`#6`#7]{{%%   %                / |
\settriparms[#1]%                            %             #5/  |
\width=\height                               %              /   |
\diagram%                                    %             v    |
\putCtrianglep<\arrowtypea`\arrowtypeb`%     %           #3     |#6
\arrowtypec;\height>%                        %             \    |
(0,0)[#2`#3`#4;#5`#6`{#7}]%                  %            #7\   |
\enddiagram%%                                %               \  |
}}                                           %                v v
\def\putDtrianglep<#1>(#2,#3)[#4`#5`#6;#7`#8`#9]{{%
\settriparms[#1]%
\xpos=#2 \ypos=#3
\advance\xpos by\height \advance\ypos by\height
\putmorphism(\xpos,\ypos)(-1,-1)[``{#9}]{\height}{\arrowtypec}r%
\advance\xpos by-\height \advance\ypos by\height
\putmorphism(\xpos,\ypos)(1,-1)[`#5`{#8}]{\height}{\arrowtypeb}r%
{\multiply\height by 2
\putvmorphism(\xpos,\ypos)[#4`#6`{#7}]{\height}{\arrowtypea}l}%
}}
 
\def\putDtriangle{\@ifnextchar <{\putDtrianglep}{\putDtrianglep
    <\arrowtypea`\arrowtypeb`\arrowtypec;\height>}}
\def\Dtriangle{\@ifnextchar <{\Dtrianglep}{\Dtrianglep
   <\arrowtypea`\arrowtypeb`\arrowtypec;\height>}}
\def\Dtrianglep<#1>[#2`#3`#4;#5`#6`#7]{{%%  %          | \
\settriparms[#1]%                           %          |  \#6
\width=\height                              %          |   \
\diagram%                                   %          |    v
\putDtrianglep<\arrowtypea`\arrowtypeb`%    %        #5|    #3
\arrowtypec;\height>%                       %          |    /
(0,0)[#2`#3`#4;#5`#6`{#7}]%                 %          |   /#7
\enddiagram%%                               %          |  /
}}                                          %          v v
\def\setrecparms[#1`#2]{\width=#1 \height=#2}%
\def\recursep<#1`#2>[#3;#4`#5`#6`#7`#8]{{\m@th
\width=#1 \height=#2
\settokens`#3`
\settowidth{\tempdimen}{$\tokena$}
\ifdim\tempdimen=0pt
  \savebox{\tempboxa}{\hbox{$\tokenb$}}%
  \savebox{\tempboxb}{\hbox{$\tokend$}}%
  \savebox{\tempboxc}{\hbox{$#6$}}%
\else
  \savebox{\tempboxa}{\hbox{$\hbox{$\tokena$}\times\hbox{$\tokenb$}$}}%
  \savebox{\tempboxb}{\hbox{$\hbox{$\tokena$}\times\hbox{$\tokend$}$}}%
  \savebox{\tempboxc}{\hbox{$\hbox{$\tokena$}\times\hbox{$#6$}$}}%
\fi
\ypos=\height
\divide\ypos by 2
\xpos=\ypos
\advance\xpos by \width
\bfig
\putCtrianglep<-1`1`1;\ypos>(0,0)[`\tokenc`;#5`#6`{#7}]%
\puthmorphism(\ypos,0)[\tokend`\usebox{\tempboxb}`{#8}]{\width}{-1}b%
\puthmorphism(\ypos,\height)[\tokenb`\usebox{\tempboxa}`{#4}]{\width}{-1}a%
\advance\ypos by \width
\putvmorphism(\ypos,\height)[``\usebox{\tempboxc}]{\height}1r%
\efig
}}
 
\def\recurse{\@ifnextchar <{\recursep}{\recursep<\width`\height>}}
 
\def\puttwohmorphisms(#1,#2)[#3`#4;#5`#6]#7#8#9{{%
\puthmorphism(#1,#2)[#3`#4`]{#7}0a
\ypos=#2
\advance\ypos by 20
\puthmorphism(#1,\ypos)[\phantom{#3}`\phantom{#4}`#5]{#7}{#8}a
\advance\ypos by -40
\puthmorphism(#1,\ypos)[\phantom{#3}`\phantom{#4}`#6]{#7}{#9}b
}}
 
\def\puttwovmorphisms(#1,#2)[#3`#4;#5`#6]#7#8#9{{%
\putvmorphism(#1,#2)[#3`#4`]{#7}0a
\xpos=#1
\advance\xpos by -20
\putvmorphism(\xpos,#2)[\phantom{#3}`\phantom{#4}`#5]{#7}{#8}l
\advance\xpos by 40
\putvmorphism(\xpos,#2)[\phantom{#3}`\phantom{#4}`#6]{#7}{#9}r
}}
 
\def\puthcoequalizer(#1)[#2`#3`#4;#5`#6`#7]#8#9{{%
\setpos(#1)%
\puttwohmorphisms(\xpos,\ypos)[#2`#3;#5`#6]{#8}11%
\advance\xpos by #8
\puthmorphism(\xpos,\ypos)[\phantom{#3}`#4`#7]{#8}1{#9}
}}
 
\def\putvcoequalizer(#1)[#2`#3`#4;#5`#6`#7]#8#9{{%
\setpos(#1)%
\puttwovmorphisms(\xpos,\ypos)[#2`#3;#5`#6]{#8}11%
\advance\ypos by -#8
\putvmorphism(\xpos,\ypos)[\phantom{#3}`#4`#7]{#8}1{#9}
}}
 
\def\putthreehmorphisms(#1)[#2`#3;#4`#5`#6]#7(#8)#9{{%
\setpos(#1) \settypes(#8)
\if a#9 %
     \vertsize{\tempcounta}{#5}%
     \vertsize{\tempcountb}{#6}%
     \ifnum \tempcounta<\tempcountb \tempcounta=\tempcountb \fi
\else
     \vertsize{\tempcounta}{#4}%
     \vertsize{\tempcountb}{#5}%
     \ifnum \tempcounta<\tempcountb \tempcounta=\tempcountb \fi
\fi
\advance \tempcounta by 60
\puthmorphism(\xpos,\ypos)[#2`#3`#5]{#7}{\arrowtypeb}{#9}
\advance\ypos by \tempcounta
\puthmorphism(\xpos,\ypos)[\phantom{#2}`\phantom{#3}`#4]{#7}{\arrowtypea}{#9}
\advance\ypos by -\tempcounta \advance\ypos by -\tempcounta
\puthmorphism(\xpos,\ypos)[\phantom{#2}`\phantom{#3}`#6]{#7}{\arrowtypec}{#9}
}}
 
\def\setarrowtoks[#1`#2`#3`#4`#5`#6]{%
\def\toka{#1}
\def\tokb{#2}
\def\tokc{#3}
\def\tokd{#4}
\def\toke{#5}
\def\tokf{#6}
}
\def\hex{\@ifnextchar <{\hexp}{\hexp<1000`400>}}
\def\hexp<#1`#2>[#3`#4`#5`#6`#7`#8;#9]{%
\setarrowtoks[#9]
\yext=#2 \advance \yext by #2
\xext=#1 \advance\xext by \yext
\bfig
\putCtriangle<-1`0`1;#2>(0,0)[`#5`;\tokb``\tokd]
\xext=#1 \yext=#2 \advance \yext by #2
\putsquare<1`0`0`1;\xext`\yext>(#2,0)[#3`#4`#7`#8;\toka```\tokf]
\advance \xext by #2
\putDtriangle<0`1`-1;#2>(\xext,0)[`#6`;`\tokc`\toke]
\efig
}
\makeatother

\newtheorem{thm}{Theorem}
\newtheorem{prop}[thm]{Proposition}
\newtheorem{cor}[thm]{Corollary}
\newtheorem{lem}[thm]{Lemma}
\newenvironment{pf}[1]{\noindent\textit{Proof.} #1}{\hfill$\Box$ \medskip}
\newenvironment{rema}[1]{\noindent {\em Remark.} #1}{}
\newenvironment{exa}[1]{\noindent {\em Example.} #1}{}

\newcommand{\ham}{\textup{Ham}}
\newcommand{\id}{\textup{id}}
\newcommand{\barespacio}{\mbox{ }|\mbox{ }}

\begin{document}

%\runningheads{A. Pedroza}{Seidel's representation on 
%$\textup{Ham}(M\times N,\omega\oplus\eta)$}

\title{
Seidel's representation on the Hamiltonian group of a Cartesian product}

\author{Andr\'es Pedroza}

\address{Facultad de Ciencias\\
         Universidad de Colima\\
     Bernal D\'{\i}az del Castillo No. 340\\
         Colima, Col., Mexico 28045}
\email{andres\_pedroza@ucol.mx}
%\corraddr{E-mail:andres\_pedroza@ucol.mx}

\begin{abstract}
Let $(M,\omega)$ be a closed symplectic manifold and 
$\textup{Ham}(M,\omega)$ the group of Hamiltonian diffeomorphisms
of $(M,\omega)$. Then the Seidel homomorphism is 
a map from the fundamental group of $\textup{Ham}(M,\omega)$
to the quantum homology ring $QH_*(M;\Lambda)$. Using this
homomorphism we give a sufficient condition for 
when a nontrivial loop $\psi$ in $\textup{Ham}(M,\omega)$ determines
a nontrivial loop $\psi\times\textup{id}_N$ in 
$\textup{Ham}(M\times N,\omega\oplus\eta)$,
where $(N,\eta)$ is a closed symplectic manifold such that $\pi_2(N)=0$.
\end{abstract}

%\thanks{The author was supported by a CONACYT grant.}
%\subjclass[2000]{Primary: 57R17 Secondary: 57S05.}

%\received{}
%\revised{}
%\accepted{}

\maketitle

\section{Introduction}

Let $(M,\omega)$ be a closed symplectic manifold and $\psi=\{\psi_t\}_{0\leq 
t\leq 1}$ a loop  about the indentity map in the group of Hamiltonian diffeomorphisms
$\ham(M,\omega)$. Then associated to $\psi$, there is a fibration
$\pi:P_\psi\to S^2$ with fiber $M$. In \cite{seidel}, P. Seidel defined
a group homomorphism  $\mathcal{S}:\pi_1(\ham(M,\omega))\to QH_*(M;\Lambda)$,
where $QH_*(M;\Lambda)$ is the quantum homology ring of $(M,\omega)$. The map $\mathcal{S}$
is usually called Seidel's representation, since its image lies in the subring of units
of $QH_*(M;\Lambda)$, which in turn defines a homomorphism of the quantum homology ring
via quantum multiplication. The homomorphism $\mathcal{S}$ can be thought as the
quantum analog of Weinstein's action $\mathcal{A}:\pi_1(\ham(M,\omega))\to \mathbb{R}/P_\omega$
of \cite{weinstein}. The element
$\mathcal{S}(\psi)$, is defined in terms of Gromov-Witten invariants related to
the moduli space of holomorphic sections of the induced fibration
$\pi:P_\psi\to S^2$. This homomorphism was used by P. Seidel to detect nontrivial loops
in the group $\ham(M,\omega)$.     % where $M$ is the complex Grassmanian. 
Further, the type of Gromov-Witten invariants involved in the definition
of $\mathcal{S}$, where studied by D. McDuff in \cite{m:s2fib}, to show
that the rational cohomology of a Hamiltonian fibration splits.

For very special symplectic manifolds the fundamental group of
$\ham(M,\omega)$ is completely known. The easiest case is when $M$ has dimension
2, since in this case symplectic diffeomorphisms agree with volume
preserving diffeomorphisms. Hence the fundamental group of
$\ham(S^2,\omega)$ is isomorphic to $\mathbb{Z}_2$; and 
 $\ham(\Sigma_g,\omega)$ is contractible for  $g\geq 1$.
For further details see \cite{pol}. In higher dimensions M. Gromov 
showed in \cite{gromov} that the fundamental group of $\ham(\mathbb{C}P^2,\omega_{FS})$
is isomorphic to $\mathbb{Z}_3$ and the fundamental group of
$\ham(S^2\times S^2,\omega\oplus\omega)$ is isomorphic to a semidirect
product of $\mathbb{Z}_2$ with itself. The last case is more interesting
when  the standard symplectic $\omega\oplus\omega$ on $S^2\times S^2$ is replaced by the 
symplectic form $\lambda\omega\oplus \omega$,
  where $\lambda$ is a real constant greater than 1.
In this case D. McDuff proved %in \cite{dusa} 
that the fundamental group
of $\ham(S^2\times S^2,\omega\oplus\lambda\omega)$ contains an element
of infinite order.

Consider $(M,\omega)$ and $(N,\eta)$ closed symplectic manifolds. Moreover
assume that both manifolds are monotone. If $\psi$ is a loop
of Hamiltonian diffeomorphisms of $(M,\omega)$ based at the identity map, then $\psi\times\textup{id}_N$
is a loop of Hamiltonian diffeomorphisms of the symplectic manifold $(M\times N,
\omega\oplus \eta)$. Then
both loops $\psi$ and $\psi\times\textup{id}_N$ induce fibrations
$\pi:P_\psi\to S^2$ and $\pi_1:P_{\psi\times\textup{id}_N}\to S^2$ with fibers
the symplectic manifolds $M$ and $M\times N$ respectively. This article 
aims on relating  Seidel's representations over 
$\pi_1(\ham(M,\omega))$ and $\pi_1(\ham(M\times N,\omega\oplus\eta))$.
We show that Seidel's representation on $\pi_1(\ham(M\times N,\omega\oplus\eta))$
restricted to elements of the form $\psi\times \textup{id}_N$ is 
essentially the same as Seidel's representation on $\pi_1(\ham(M,\omega))$.

We achieve this by relating the fundamental groups of $\ham(M,\omega)$
and $\ham(M\times N,\omega\oplus\eta)$   and the quantum homology rings of 
$M$ and $M\times N$. For, let 
$$
\tau:\pi_1(\ham(M,\omega))\to \pi_1(\ham(M\times N,\omega\oplus\eta))
$$ 
be the group homomorphism defined by $\tau(\psi)=\psi\times\textup{id}_N$, 
and 
$$
\kappa:QH_*(M;\Lambda)\to QH_{*+\textup{dim}(N)}(M\times N;\Lambda)
$$
be the map defined  on homogeneous elements %  $\alpha\otimes q^st^r$
by $\kappa(\alpha\otimes q^st^r)=(\alpha\otimes[N])\otimes q^st^r$
where $\alpha\in H_*(M)$ and  $\alpha\otimes[N]\in H_*(M\times N)$.  Extend  $\Lambda$-linearly
the map $\kappa$ to all the quantum homology ring $QH_*(M;\Lambda)$. 
In Section \ref{s:quantkun} we 
show that $\kappa$ is in fact a ring homomorphism under the quantum product,
if both manifolds $(M,\omega)$ and $(N,\eta)$ are monotone with the same
constant. This statement, as others in this article, is a direct consequence of the fact 
that the quantum homology ring $QH_*(M\times N;\Lambda)$ satisfies the K\"unneth formula.

\begin{thm}
\label{t:commdiag}
Let $(M,\omega)$ and $(N,\eta)$ be closed symplectic manifolds. Assume that
$(M,\omega)$ has dimension $2n$ and is monotone,  and that $\pi_2(N)=0$. Then 
 $\mathcal{S}\circ\tau(\psi)=\kappa\circ\mathcal{S}(\psi)$
for every
$\psi$ in $\pi_1(\ham(M,\omega))$.
That is the following diagram
commutes

\begin{center}
\begin{picture}(2500,700)(0,0)
\put(0,500){$\pi_1(\ham(M,\omega)) \hskip 3cm
  \pi_1(\ham(M\times N,\omega\oplus\eta))$}
\put(70,10){$QH_{2n}(M;\Lambda) \hskip 4cm
  QH_{2n+\textup{dim}(N)}(M\times N;\Lambda)$}
\put(700,40){\vector(1,0){800}}
\put(790,520){\vector(1,0){540}}
\put(350,430){\vector(0,-1){260}}
\put(2100,430){\vector(0,-1){260}}
\put(190,250){\begin{small}$\mathcal{S}$ \end{small}}
\put(2180,250){\begin{small}$\mathcal{S}$ \end{small}}
\put(1000,600){\begin{small}$\tau$ \end{small}}
\put(1000,120){\begin{small}$\kappa$ \end{small}}
\end{picture}
\end{center}
\end{thm}

It is important to  relate Thm. \ref{t:commdiag} with a result of D. McDuff and S. Tolman. 
In \cite{mt}, McDuff and Tolman found a formula for $\mathcal{S}(\psi)$ in the case
when $\psi$ is a Hamiltonian circle action on $(M,\omega)$. Thus if $\psi$ is a Hamiltonian 
circle action on $(M,\omega)$,  
denote by $K:M\to \mathbb{R}$ the normalized moment map of the circle action and by
$K_0$  the maximum value of $K$. Let  $M_\textup{max}$ be the symplectic 
 submanifold on which the moment map $K$ 
achieves its maximun. Note that $M_\textup{max}$ is part of the fixed point set of the
circle action.  Finally assume that there
is a neighborhood $U$ of $M_\textup{max}$ such that the action of the circle 
is free on $U-M_\textup{max}$. Under the above assumptions  McDuff-Tolman formula
for the circle action $\psi$ reads
\begin{eqnarray}
\label{fmt}
\mathcal{S}(\psi)&=&[M_\textup{max}]\otimes q^{\textup{codim}(M_\textup{max})/2}
t^{-K_0} + \\
& &
\sum_{ \{A:  \tilde\omega(A)> K_0 \} }  \alpha_{A}\otimes q^{-c(A)} t^{-\tilde\omega(A)}.
\nonumber
\end{eqnarray}
In order to make a clear statement of the goal of this article we postpone to Section \ref{s:smallqh}, 
the definitions of the elements $A,c(A)$ and $\tilde\omega(A)$ that appear in Eq. (\ref{fmt}).
% and also those that appear in Eq. (\ref{f2mt}).
Also a word of warning about Eq. (\ref{fmt}). We have used the notation of \cite{msjholo}
in stating McDuff-Tolman formula and not that of the original paper \cite{mt}.
At then end of Section \ref{s:smallqh}, we clarify how they are related. 

Now let $(M,\omega)$ and $\psi$ be as above and $(N,\eta)$ any closed symplectic 
manifold. Then $\psi\times \textup{id}_N$
is also a Hamiltonian circle action on $M\times N$ with moment map $H:M\times N\to \mathbb{R}$
given by $H(p,q)=K(p)$. So defined, the moment map $H$  is normalized. 
Also $(M\times N)_\textup{max}=M_\textup{max}\times N$
and $H_0=K_0$. 
Observe that $\textup{codim}_M(M_\textup{max})=\textup{codim}_{M\times N}(M_\textup{max}\times N)$.
Thus McDuff-Tolman formula for the circle action $\psi\times \textup{id}_N$ is,
\begin{eqnarray}
\label{f2mt}
\mathcal{S}(\psi\times \textup{id}_N)&=&[M_\textup{max}\times N]\otimes q^{\textup{codim}(M_\textup{max})/2}
t^{-K_0} +\\
& & \sum_{\{A^\prime:\tilde\omega^\prime(A^\prime)> K_0\}} \alpha^\prime_{A^\prime}\otimes q^{-c^\prime(A^\prime)} 
t^{-\tilde\omega^\prime(A^\prime)} \nonumber.
\end{eqnarray}

At this point is important to observe that the first term on the right hand side 
of Eqs. (\ref{fmt}) and (\ref{f2mt}) only differ by the class $[N]$.
Well  Thm. \ref{t:commdiag} guarantees that if $\pi_2(N)=0$ not only the
first terms of $\mathcal{S}(\psi)$ and $\mathcal{S}(\psi\times \textup{id}_N)$
differ by $[N]$, but the equality  $\mathcal{S}(\psi)
\otimes [N]= \mathcal{S}(\psi\times \textup{id}_N)$ holds  in  $QH_*(M\times N;\Lambda)$
for any loop $\psi$, not just a circle action.

\bigskip

Notice that the map $\kappa:QH_*(M;\Lambda)\to QH_{*+\textup{dim}(N)}(M\times N;\Lambda)$
so defined is injective. Therefore Thm. \ref{t:commdiag} tells us
when a nontrivial $\psi \in \pi_1(\ham(M,\omega))$ induces
a nontrivial element $\psi\times \textup{id}_N$ in $\pi_1(\ham(M\times N,\omega\oplus\eta))$.

\begin{cor}
\label{c:senprod}
Let $(M,\omega)$ and $(N,\eta)$ be as in Thm. \ref{t:commdiag}. Then if
$\psi\in \pi_1(\ham(M,\omega))$ is such that $\mathcal{S}(\psi)\neq 1=[M]$, then
the loop
$\psi\times \textup{id}_N$ is also nontrivial in $\pi_1(\ham(M\times N,\omega\oplus\eta))$.
\end{cor}

Hence if the Seidel representation on $\pi_1(\ham(M,\omega))$ is injective,
%like in the case of $S^2$, 
we conclude that the group homomorphism  
$\tau:\pi_1(\ham(M,\omega))\to \pi_1(\ham(M\times N,\omega\oplus\eta))$ 
is also injective. For instance, by the result of McDuff and Tolman we
know that the Seidel representation is injective in the case  when $M=S^2$
or $\mathbb{C}P^2$.
Therefore for any closed symplectic manifold $(N,\eta)$ such that
$\pi_2(N)=0$, we have that the group homomorphism
$\tau:\pi_1(\ham(M,\omega))\to \pi_1(\ham(M\times N,\omega\oplus\eta))$ 
is injective for $M=S^2$ or $\mathbb{C}P^2$.

\medskip
\begin{exa}
Let $(M,\omega)$ and $(N,\eta)$ be symplectic manifolds as in Thm. \ref{t:commdiag}.
Moreover assume that there is a loop $\gamma$ in $\textup{Ham}(M,\omega)$ such 
that $\mathcal{S}(\gamma)$ has infinite order in $QH_*(M;\Lambda)$ under quantum 
multiplication. Thus the loop $\gamma$
also has infinite order in the fundamental group of $\textup{Ham}(M,\omega)$.

Hence we have that $\mathcal{S}(\gamma^m)$ is not equal to the identity 
$1=[M]$ in $QH_*(M;\Lambda)$
for all $m\in \mathbb{Z}$ different from zero. Then by Cor. \ref{c:senprod},
it follows that the loop $\gamma^m\times \textup{id}_N$ is not homologous to
the constant loop in $\ham(M\times N,\omega\oplus\eta)$ for all nonzero $m$.
That is, $\gamma^m\times \textup{id}_N$ is  an element of infinite order in
$\pi_1(\ham(M\times N,\omega\oplus\eta))$.
\end{exa}

\bigskip 
Now consider the case when $M=N$. Hence assume that $(M,\omega)$ is a closed
symplectic manifold such that $\pi_2(M)$ is trivial. Thus $(M,\omega)$
is monotone. We are intrested in understanding when
a nontrivial loop $\psi$ in $\textup{Ham}(M,\omega)$ induces a nontrivial loop
$\psi\times \psi$ in $\textup{Ham}(M\times M,\omega\oplus\omega)$. 
That is, we are intrested in the image of the group homomorphism
$$
\tau^\prime:\pi_1(\ham(M,\omega))\to \pi_1(\ham(M\times M,\omega\oplus\omega))
$$ 
defined by $\tau^\prime(\psi)=\psi\times\psi$.

Consider the map
$$
\kappa^\prime:QH_*(M;\Lambda)\to QH_{*+*}(M\times M;\Lambda)\simeq (H_{*}(M)
\otimes_\mathbb{Z} H_{*}(M))\otimes_\mathbb{Z}\Lambda
$$
%by $\kappa^\prime(x)=x\otimes x$ for $x\in QH_*(M;\Lambda)$. That is 
defined on homogeneous elements by $\kappa^\prime(\alpha\otimes q^st^r)=(\alpha\otimes \alpha)
\otimes q^{2s}t^{2r}$ where $\alpha\in H_*(M)$.
%and extend it additively to all $QH_*(M;\Lambda)$. 
In Section \ref{s:quantkun} we will review the fact that the K\"unneth
formula holds in quantum homology. Hence the map $\kappa^\prime$ corresponds
to the diagonal map
$$
\Delta:QH_*(M;\Lambda)\to  QH_*(M;\Lambda)\otimes_\Lambda QH_*(M;\Lambda)\simeq
QH_{*}(M\times M;\Lambda)
$$
defined by $\Delta(x)=x\otimes x$ for all $x\in QH_*(M;\Lambda)$
via the quantum K\"unneth formula.

\begin{thm}
\label{t:2}
Let $(M,\omega)$ be closed symplectic manifold of dimension $2n$ such that $\pi_2(M)$ is trivial.
Then
 $\mathcal{S}\circ\tau^\prime(\psi)=\kappa^\prime\circ\mathcal{S}(\psi)$
 for every
$\psi$ in $\pi_1(\ham(M,\omega))$.
That is the following diagram
commutes

\begin{center}
\begin{picture}(2500,700)(0,0)
\put(0,500){$\pi_1(\ham(M,\omega)) \hskip 3cm
  \pi_1(\ham(M\times M,\omega\oplus\omega))$}
\put(70,10){$QH_{2n}(M;\Lambda) \hskip 4cm
  QH_{4n}(M\times M;\Lambda)$}
\put(700,40){\vector(1,0){800}}
\put(790,520){\vector(1,0){540}}
\put(350,430){\vector(0,-1){260}}
\put(2100,430){\vector(0,-1){260}}
\put(190,250){\begin{small}$\mathcal{S}$ \end{small}}
\put(2180,250){\begin{small}$\mathcal{S}$ \end{small}}
\put(1000,600){\begin{small}$\tau^\prime$ \end{small}}
\put(1000,120){\begin{small}$\kappa^\prime$ \end{small}}
\end{picture}
\end{center}
\end{thm}

As in the case of Thm. \ref{t:commdiag}, one can use McDuff-Tolman formula to verify that
 Thm. \ref{t:2} works in the case when the Hamiltonian loop $\psi$ is a Hamiltonian circle
action. 
In fact,  if $(M,\omega)$, $\psi$, and $K:M\to \mathbb{R}$ are as above, that is $\psi$ is a
Hamiltonian circle action with normalized moment map $K$, then  $\psi\times \psi$
is a Hamiltonian circle action on the product manifold $(M\times M,\omega\oplus
\omega)$ with normalized moment map $H:M\times M\to \mathbb{R}$ given by 
$H(p_1,p_2)=K(p_1)+K(p_2)$. Hence the maximun value $H_0$ of the
moment map $H$ satisfies the relation  $H_0=2K_0$  
and also $(M\times M)_\textup{max}=M_\textup{max}\times M_\textup{max}$.
Then in this case McDuff-Tolman formula for the Hamiltonian circle
action $\psi\times\psi$ on $M\times M$ is given by
\begin{eqnarray}
\label{f3mt}
\mathcal{S}(\psi\times\psi)&=&[M_\textup{max}\times M_\textup{max}]\otimes q^{\textup{codim}_{M\times M}(
M_\textup{max}\times M_\textup{max})/2}
t^{-H_0} + \nonumber\\
& &
\sum_{ \{A:  \tilde\omega(A)> H_0 \} }  \alpha_{A}\otimes q^{-c(A)} t^{-\tilde\omega(A)}\nonumber\\
&=&[M_\textup{max}\times M_\textup{max}]\otimes q^{\textup{codim}_{M}(
M_\textup{max})}
t^{-2K_0} + \nonumber\\
& &
\sum_{ \{A:  \tilde\omega(A)> 2K_0 \} }  \alpha_{A}\otimes q^{-c(A)} t^{-\tilde\omega(A)}.
\end{eqnarray}

Comparing the first term on the right hand side of Eqs. (\ref{fmt}) and (\ref{f3mt}), one 
checks that they are related by the map $\kappa^\prime$. Well according to Thm. \ref{t:2},
$\mathcal{S}(\psi)$ and $\mathcal{S}(\psi\times\psi)$ are related by the map
$\kappa^\prime$ for any loop $\psi$ in $\textup{Ham}(M,\omega)$, not just  a
Hamiltonian  circle action.

As before the map $\kappa^\prime$ so defined is injective. Hence we have a criteria to
determine when the loop $\psi\times \psi$ is nontrivial in $\pi_1(\ham(M\times M,\omega\oplus\omega))$.

\begin{cor}
 Let $(M,\omega)$ be a closed symplectic manifold such that 
 $\pi_2(M)=0$. If 
$\psi\in \pi_1(\ham(M,\omega))$ is such that $\mathcal{S}(\psi)\neq 1=[M]$, then
the Hamiltonian loop $\psi\times \psi$ is nontrivial in 
$\pi_1(\ham(M\times M,\omega\oplus\omega))$.
\end{cor}

\bigskip

The author would like to thank Prof. Dusa McDuff  for her patience on reading the first draft of the manuscript and making valuable observations to it; and the Referee for the useful comments and suggestions.
The author was partially supported  by 
Consejo Nacional de Ciencia y Tecnolog\'ia M\'exico grant.

%%%%%%%%%%%%%%%%%%%%%%%%%%%%%%%%%%%%%%%%%%%%%%%%
\section{Hamiltonian fibrations}
%%%%%%%%%%%%%%%%%%%%%%%%%%%%%%%%%%%%%%%%%%%%%%%%

Consider $(M,\omega)$ a closed symplectic manifold. Let $\psi=\{\psi_t\}_{0\leq t\leq 1}$
 be a loop  about the indentity map in the group of
Hamiltonian diffeomorphisms $\ham(M,\omega)$. Associated to $\psi$ there is a 
smooth fibration $\pi:P_\psi\to S^2$ with fiber $M$ defined as follows.
Let $D^+=\{z\in \mathbb{C} \barespacio | z|\leq 1\}$ be the closed unit disk with
the positive orientation. Then the total space of the fibration  $P_\psi$ is defined as
$$
(D^+\times M) \amalg (D^-\times M) / \sim
$$
where $(e^{it},p)_+\sim (e^{-it},\psi_t^{-1}(p))_-$ and $D^-$ stands for $D^+$ with the
opposite orientation. This fibration has $\ham(M,\omega)$ as structure group, 
and is called  \textbf{Hamiltonian fibration}. See (\cite{msjholo}, p. 251).

In fact there is a one-to-one correspondence between homotopic loops in 
$\ham(M,\omega)$ based at the identity map $\textup{id}_M$ and isomorphic
fibrations over $S^2$ with fiber $M$ and structure group $\ham(M,\omega)$. 
In order to avoid cumbersome notation
we will use the same notation to denote  loops based at the identity in
$\ham(M,\omega)$ and its homotopy class in $\pi_1(\ham(M,\omega))$, namely 
$\psi=\{\psi_t\}_{0\leq t\leq 1}$.

In a Hamiltonian fibration $\pi:P_\psi\to S^2$ with fiber the symplectic
manifold $(M,\omega)$, there exists a 
closed $2$-form $\tilde\omega$ on $P_\psi$ such that it restricts to  $\omega_z$ on every fiber
$(P_\psi)_z$  for $z\in S^2$, and such that 
%Moreover one normalize the  form $\tilde\omega$ by requiring that
$$
\pi_!(\tilde\omega^{n+1})=0,
$$
where $\pi_!: H^*(P_\psi)\to H^{*-\textup{dim(M)}}(S^2)$ stands for integration along the fiber $M$. A 
 $2$-form $\tilde\omega$ that satisfies the above conditions
is called a \textbf{coupling form}. See \cite{gs} for more details. 
The coupling form  $\tilde\omega$
defines a connection on the fibration, where
the horizontal distribution is defined as the $\tilde\omega$-complement
of the vertical subspace. That is, for $p\in P_\psi$,
$$
\textup{hor}(T_pP_\psi)=\{v\in T_pP_\psi\barespacio \tilde\omega(v,u)=0 \textup{ for all }
u\in \textup{ker} (\pi_{*,p})\}.
$$

There is another canonical class associated to the fibration $\pi:P_\psi\to S^2$,
apart from the cohomology class determined by
coupling form. Recall that the vertical vector bundle of a fibration is the vector bundle
$$
T^\textup{V}P_\psi=\{v \in T_pP_\psi \barespacio\pi_{*,p}(v)=0\}
$$
over the total space $P_\psi$. The coupling form $\tilde\omega$
restricted to this subbundle is nondegenerate. Thus the first Chern
class of $T^\textup{V}P_\psi$ is well defined. This class is denoted 
  by $c_\psi:=c_1(T^\textup{V}P_\psi)   \in H^2(P_\psi;\mathbb{Z})$.

\medskip

Let $(M,\omega)$ and $(N,\eta)$ be closed symplectic manifolds,
then $(M\times N,\omega\oplus\eta)$ is also a closed symplectic
manifold. 
Note that the $2$-form $\omega\oplus\eta$ is a shorthand notation
for the $2$-form $(\textup{pr}_M)^*(\omega)+(\textup{pr}_N)^*(\eta)$  
on $M\times N$, 
where $\textup{pr}_M$
and $\textup{pr}_N$ are the projection maps from $M\times N$ to $M$ and
$N$ respectively.
If  $\psi$  is a loop in the group $\ham(M,\omega)$, then
$\psi\times \id_N$ is also a loop of Hamiltonian diffeomorphisms of
the symplectic manifold $(M\times N,\omega\oplus\eta)$.
Thus there
is a Hamiltonian fibration $\pi_1:P_{\psi\times \textup{id}_N}\to
S^2$ with fiber $M\times N$. As before we have the fiber bundle
 $\pi:P_\psi\to S^2$ with fiber $M$.
% but from now on we denote the map by $\pi_1$ instead of $\pi.$
 Define the fibration
$\pi_0:P_\psi\times N\to S^2$ where the projection map is defined as $\pi_0(x,q)=\pi(x)$. So defined
$\pi_0$ is a fiber bundle with fiber $M\times N$. Well both fiber
bundles $P_{\psi\times \textup{id}_N}$ and $P_\psi\times N$ are isomorphic.

\begin{lem} \label{l:isofibrations}
The fiber bundles $\pi_0:P_\psi\times N\to S^2$ and  $\pi_1:P_{\psi\times \textup{id}_N}
\to S^2$ with fiber $M\times N$ are isomorphic fibrations. 
%That is there are fibrewise diffeomorphic.
\end{lem}
\begin{pf}
Consider the map
$
\rho:
\left( (D^+\times M) \amalg (D^-\times M) \right)\times N \to  P_{\psi\times \textup{id}_N}
%(D^+\times M\times N) \amalg (D^-\times M\times N)
$
defined as $\rho((u,p),q)=[u,p,q]$. Then
$$
\rho((e^{it},p),q)=[e^{it},p,q]=[e^{-it},\psi_t^{-1}(p),\textup{id}_N(q)]
=
\rho((e^{-it},\psi_t^{-1}(p)),q).
$$
This means that $\rho$ induces a map on the quotient $P_\psi\times N$.
We denote such map by the same letter $\rho$. So defined
$
\rho:
P_\psi\times N \to  P_{\psi\times \textup{id}_N}
$
is smooth and bijective. Moreover for any $([u,p],q)\in P_\psi\times N$
we have that
$$
\pi_0 ([u,p],q)=\pi([u,p])=u 
\hbox{ \hskip .6cm and \hskip .6cm}
\pi_1\circ\rho ([u,p],q)=
\pi_1 ([u,p,q])= u.
$$
Therefore  $\pi_0=\pi_1\circ\rho$ and $\rho$ is fiberwise preserving. That is the fiber
bundles $P_\psi\times N$ and $P_{\psi\times \textup{id}_N}$ are isomorphic fibrations.
\end{pf}

Thus there are two isomorphic fibrations over $S^2$ with fiber $M\times N$;
 $P_{\psi\times \textup{id}_N}$ and $P_\psi\times N$. Hence both  vertical
bundles are isomorphic $T^\textup{V}(P_{\psi\times \textup{id}_N})\simeq 
T^\textup{V}(P_\psi\times N)$ as vector bundles. In order to compare the
first Chern classes of both fibrations,  consider the projections 
maps $\lambda_1: P_\psi\times N \to P_\psi$, $\lambda_2: P_\psi\times N \to N$ 
and the diagram

\begin{center}
\begin{picture}(2200,700)(0,0)
\put(-30,500){$T^\textup{V}P_{\psi}$}
\put(10,10){$P_\psi$}
\put(600,500){ $(\lambda_1)^{*}(T^\textup{V}P_{\psi})\oplus (\lambda_2)^{*} (TN)$ }
\put(1000,10){$P_\psi\times N$}
\put(2240,500){$TN$}
\put(2240,10){$N$}
\put(940,40){\vector(-1,0){700}}
\put(1400,40){\vector(1,0){700}}
\put(540,540){\vector(-1,0){300}}
\put(1850,540){\vector(1,0){300}}
\put(50,430){\vector(0,-1){260}}
\put(1180,430){\vector(0,-1){260}}
\put(2300,430){\vector(0,-1){260}}
\put(550,70){\begin{small}$\lambda_1$ \end{small}}
\put(1790,70){\begin{small}$\lambda_2$ \end{small}}
\end{picture}
\end{center}
\medskip
where $(\lambda_1)^{*}(T^\textup{V}P_{\psi})$ and $(\lambda_2)^{*} (TN)$ stand for the
pullback bundles.

\begin{prop}
\label{p:vertiso}
The  vector bundles  %$(P_{\psi\times \textup{id}_N})^\textup{v}\simeq 
$T^\textup{V}(P_\psi\times N)$ and
$(\lambda_1)^{*}(T^\textup{V} P_{\psi})\oplus (\lambda_2)^{*} (TN)$ 
 %over $P_\psi\times N\simeq P_{\psi\times \textup{id}_N}$ 
are  isomorphic vector bundles over $P_\psi\times N$.
\end{prop}
\begin{pf}
Let $x=(p,q)\in P_\psi\times N$ and $u=u_1+u_2\in T_x(P_\psi\times N)\simeq 
T_pP_\psi\oplus T_qN$. Thus the vector $u$ belongs to  $T^\textup{V}(P_\psi\times N)$
if and only if $(\pi_0)_*(u)=0$.  By the definition of the map $\pi_0$, we have
 $(\pi_0)_{*,x}(u_1+u_2)=(\pi)_{*,p}(u_1)=0$.
Thus $u=u_1+u_2$ is in $T^\textup{V}(P_\psi\times N)$ if and only if
$u_1$ is in $T^\textup{V}P_\psi$ and $u_2$ is in $TN$. Therefore
$$T^\textup{V}(P_\psi\times N)\simeq
(\lambda_1)^{*}(T^\textup{V}P_{\psi})\oplus (\lambda_2)^{*} (TN)
$$ 
as vector bundles.
\end{pf}

Let $c_\psi\in H^2(P_\psi;\mathbb{Z})$ be the first Chern class of the vector bundle 
$T^\textup{V}P_\psi\to P_\psi$. And respectively
$c_{\psi\times \textup{id}_N}\in H^2(P_{\psi\times\textup{id}_N};\mathbb{Z})
=H^2(P_{\psi}\times N;\mathbb{Z})$.

\begin{lem} 
\label{l:samec}
On $H^2(P_{\psi\times \textup{id}_N};\mathbb{Z})$ we have the identity
$$
c_{\psi\times \textup{id}_N}=(\lambda_1)^*(c_\psi)+ (\lambda_2)^*(c_1(N))
$$
where $c_1(N)$ stands for the first Chern class of $(N,\eta)$.
\end{lem}
\begin{pf}
By definition we have
$c_{\psi\times \textup{id}_N}=  c_1(T^\textup{V}P_{\psi\times \textup{id}_N})$.
Then it follows by  Prop. \ref{p:vertiso} that
\begin{eqnarray*}
c_{\psi\times \textup{id}_N}&=&  c_1(T^\textup{V}P_{\psi\times \textup{id}_N}) \\
	&=&c_1( (\lambda_1)^{*}(T^\textup{V}P_{\psi})\oplus (\lambda_2)^{*} (TN)) \\
        &=&c_1( (\lambda_1)^{*}(T^\textup{V}P_{\psi}) )+ c_1((\lambda_2)^{*} (TN)) \\
         &=& (\lambda_1)^{*} (c_\psi) + (\lambda_2)^{*} (c_1(N)).
\end{eqnarray*}
\end{pf}

A similar result holds for the coupling forms of the fibrations
$P_\psi$ and $P_{\psi\times \textup{id}_N}$.

 \begin{lem}
\label{l:samecoup}
If  $\tilde\omega$ is a coupling form of the Hamiltonian fibration
$\pi:P_\psi\to S^2$, then 
$(\lambda_1)^*(\tilde\omega)+ (\lambda_2)^*(\eta)$
is a coupling 
form of the fibration $\pi_0: P_{\psi}\times N\to S^2$.
 \end{lem}
The proof follows from Lemma \ref{l:isofibrations} and the definition
of the coupling form. We will write $\tilde\omega\oplus\eta$ for
the coupling form $(\lambda_1)^*(\tilde\omega)+ (\lambda_2)^*(\eta)$
on $P_\psi\times N\simeq P_{\psi\times\textup{id}_N}.$

\medskip
\begin{rema} 
In the proof of the main theorem it will be important to note the following.
Let  $A\in H_2(P_{\psi\times \textup{id}_N}; \mathbb{Z})$ be any class 
such that $(\lambda_2)_*(A)=0$. Then it follows from Lemmas \ref{l:samec}
and  \ref{l:samecoup} that
\begin{eqnarray}
 \label{eq.samec}
c_{\psi\times \textup{id}_N}(A)=
c_\psi((\lambda_1)_*(A)),
\end{eqnarray}
and 
\begin{eqnarray}
 \label{eq.samew}
\tilde\omega\oplus\eta(A)=
\tilde\omega((\lambda_1)_*(A)).
\end{eqnarray}
\end{rema}

%%%%%%%%%%%%%%%%%%%%%%%%%%%%%%%%%%%%%%%%%%%%%%%%%%%%%%%%%%%%%%%
%%%%%%%%%%%%%%%%%%%%%%%%%%%%%%%%%%%%%%%%%%%%%%%%%%%%%%%%%%%%%%%
\section{Small quantum homology and Seidel's homomorphism}
%%%%%%%%%%%%%%%%%%%%%%%%%%%%%%%%%%%%%%%%%%%%%%%%%%%%%%%%%%%%%%%
%%%%%%%%%%%%%%%%%%%%%%%%%%%%%%%%%%%%%%%%%%%%%%%%%%%%%%%%%%%%%%%
\label{s:smallqh}

In this section we will review the concepts
 needed to define Seidel's representation. We will follow
closely the exposition and notations of D. McDuff and  D.  Salamon
\cite{msjholo}.

Let $\psi$ be a loop in the group of Hamiltonian diffeomorphisms
of $(M,\omega)$ and $\pi:P_\psi\to S^2$  the Hamiltonian fibration associated  to
the loop $\psi$. Consider $\tilde{\omega}$ a coupling form on the fibration. 
Then  for a large positive constant $K$ the form $\Omega:= \tilde{\omega}
+ K\pi^*(\omega_0)$
on $P_\psi$  is  a nondegenerate $2$-form, where  $\omega_0$ is an area form on $S^2$. 
Is important to note that $\Omega$ and $\tilde \omega$ induced the
same horizontal distribution on $P_\psi$. Denote by 
$\mathcal{J}(P_\psi,\pi,\Omega)$ the set of almost complex
structures $J$ on $P_\psi$ that are $\Omega$-compatible and such 
that the projection map $\pi:(P_\psi,J)\to (S^2,j_0)$ is
holomorphic. Here $j_0$ is an arbitrary complex structure on $S^2$.
 Recall that $\Omega$-compatible means in particular that
$\Omega(Ju,Jv)=\Omega(u,v)$; and $\pi:(P_\psi,J)\to (S^2,j_0)$ is
holomorphic if $j_0\circ(d\pi)=(d\pi)\circ J$.
Since for any $J\in\mathcal{J}(P_\psi,\pi,\Omega)$ the projection map $\pi$ is holomorphic, 
then $J$ preserves the vertical tangent space of $P_\psi$. Also since
$J$ is a $\Omega$-compatible almost complex structure, $J$ preserves the horizontal 
distribution of $P_\psi$.

Consider a spherical class $A\in H_2(P_\psi;\mathbb{Z})$, that is
$A$ is in the image of the Hurewicz homomorphism $\pi_2(P_\psi)\to H_2(P_\psi)$.
Then if $J\in
\mathcal{J}(P_\psi,\pi,\Omega)$, let
$\mathcal{M}(A; J)$ be the moduli space of $J$-holomorphic sections of $P_\psi$ 
that represent the class $A$,
$$
\mathcal{M}(A; J) =\{ u:S^2\to P_\psi \barespacio \bar{\partial}_{J}(u)=0, \pi\circ
u=\textup{id}_{S^2}, [u]=A \}.
$$
Here $\bar{\partial}_{J}$ stands for the Cauchy-Riemann equation
$$
\bar{\partial}_{J}(u)=\frac{1}{2}\left(
du+J\circ du\circ j_0\right).
$$
where $j_0$ is a fix complex structure on $S^2$.

To assure that the moduli space $\mathcal{M}(A;  J)$ is a smooth finite dimensional 
manifold, one considers the linearized operator 
$$
D_{u}: \Omega^0(S^2,u^*(TP_\psi)) \to \Omega^{0,1}(S^2,u^*(TP_\psi))
$$
of $\bar{\partial}_{J}$.
One finds that there is a subset $\mathcal{J}_{reg}(P_\psi,\pi,\Omega)
\subset\mathcal{J}(P_\psi,\pi,\Omega)$ such that if $J$ is in
$\mathcal{J}_{reg}(P_\psi,\pi,\Omega)$, then $\mathcal{M}(A; J)$
is a smooth manifold of dimension  
$2n+2c_\psi(A),$ where the symplectic manifold $(M,\omega)$ has  dimension 
$2n$. Moreover $\mathcal{M}(A; J)$ carries a natural orientation.
The set $\mathcal{J}_{reg}(P_\psi,\pi,\Omega)$ is characterized by the fact
that  $J$ is regular if and only if for every holomorphic $J$-curve $u\in \mathcal{M}(A; J)$ the operator
$D_u$ is surjective.
For the details see (\cite{msjholo}, Ch. 3).

Denote by $\mathcal{M}_k(A; J)$ be the moduli space of $J$-holomorphic sections with
$k$ marked points. That is
$$
\mathcal{M}_k(A;J)=\{(u,z_1,\ldots,z_k)\barespacio u\in\mathcal{M}(A; J), z_i\in S^2, 
z_i\neq z_j, \forall i\neq j \}.
$$
This moduli space has dimension $\mu:=2n+2c_\psi(A)+2k$. Consider the evaluation map
$
ev:\mathcal{M}_k(A;J)\to (P_\psi)^k
$
given by $ev(u,z_1,\ldots,z_k)=(u(z_1),\ldots,u(z_k))$.  
Now one would like the map  $ev$ to represent a cycle in
the homology of $(P_\psi)^k$. Actually the map $ev$
is a pseudocycle if the manifold $P_\psi$ is monotone. In the case at hand, 
Hamiltonian fibrations, is enough to impose this
condition on the fiber $M$ rather than on the whole manifold
$P_\psi$. A symplectic manifold $(M,\omega)$ is said to be 
\textbf{monotone} if there is $\lambda>0$ such that
$$
\omega(A)=\lambda c_1(A)
$$
for all $A\in \pi_2(M)$. Then
 if $(M,\omega)$ is monotone,   if follows that $ev$ is a 
pseudocycle of dimension $\mu$ in $(P_\psi)^k$. That is, it defines a homology class
in $H_*((P_\psi)^k)$ of degree $\mu$. 

With this at hand we can define  the corresponding Gromov-Witten invariants
of $P_\psi$. However in order to define Seidel's representation one studies  
the moduli space of sections with one \textit{fixed} marked point. 
Fix a point $z_0$ in the base $S^2$ and let $\iota:M\to P_\psi$ be the inclusion
of the fiber above the base point $z_0$. Then  the moduli space
$$
\mathcal{M}_1^\textbf{w}(A;J)=\{(u,z_0)\barespacio u\in\mathcal{M}(A; J) \}.
$$ 
is a smooth oriented manifold of dimension $2n+2c_\psi(A)$,
where $\textbf{w}=\{z_0\}$ stands for the fixed marked point. Moreover the evaluation
map
$
ev_\textbf{w}:\mathcal{M}^\textbf{w}_1(A;J)\to P_\psi
$
is a pseudocycle in $P_\psi$. If we consider the inclusion map $\iota:M\to P_\psi$, then
$\iota^{-1}\circ ev_\textbf{w}$ represents a pseudocycle in $M$ of degree $2n+2c_\psi(A)$.

Let $H_*(M)$ be the torsion-free part of the group $H_*(M;\mathbb{Z})$.
Then the Gromov-Witten invariant is defined as the homomorphism
$\textup{GW}^{P_\psi,\textbf{w}}_{A,1}: H_{-2c_\psi(A)}(M) \to \mathbb{Z}$ 
given by
$$
\textup{GW}^{P_\psi,\textbf{w}}_{A,1}(\alpha)= (\iota^{-1}\circ ev_\textbf{w})\cdot_M\alpha
$$
where $\cdot_M$ denotes the cycle intersection product in $H_*(M)$ and
$A\in H_2(P_\psi;\mathbb{Z})$ a spherical class. Geometrically,
 $\textup{GW}^{P_\psi,\textbf{w}}_{A,1}(\alpha)$ is the number
of holomorphic sections of $\pi: P_\psi\to S^2$  such that $u(z_0)$
lies in the cycle $X$,  where $\alpha=[X]\in H_*(M)$.

Consider the universal Novikov ring $\Lambda^{\textup{univ}}$ defined as
$$
\Lambda^{\textup{univ}}=\left\{\sum_{s\in \mathbb{R}} r_st^s \barespacio r_s\in\mathbb{Z},
\#\{s>c: r_s\neq 0 \} <\infty \mbox{ for every } c\in \mathbb{R}
 \right\}
$$
and the graded polynomial ring $\Lambda:=\Lambda^\textup{univ}[q,q^{-1}]$
 where $q$ has degree 2. Then the small
quantum homology of $(M,\omega)$ with coefficients in $\Lambda$ is 
defined as
$$
QH_*(M;\Lambda):=H_*(M)\otimes_\mathbb{Z} \Lambda.
$$
Actually, $QH_*(M;\Lambda)$ is a ring under quantum product. The ring structure
of $QH_*(M;\Lambda)$ will be describe below.

Then Seidel's homomorphism 
$
\mathcal{S}:\pi_1(\ham(M,\omega))\to QH_*(M;\Lambda)
$
is defined as
\begin{eqnarray}
\label{e:seidel}
 \mathcal{S}(\psi)=\sum_{A\in H^{sec}_2(P_{\psi})} \mathcal{S}_A(\psi)\otimes q^{-c_\psi(A)} 
t^{-\tilde\omega_\psi(A)}
\end{eqnarray}
where the sum runs over all  spherical classes $A$ that can be realized by a section. That is,
a section $u:S^2\to P_\psi$ such that $[u]=A$.  And  $\mathcal{S}_A(\psi)$ is the
homology class in $H_{2n+2c_\psi(A)}(M)$ 
determined  by the relation
$$
\textup{GW}^{P_\psi,\textbf{w}}_{A,1}(\alpha)=\mathcal{S}_A(\psi)\cdot_M \alpha
$$
for all $\alpha\in H_*(M).$  Observe that $\mathcal{S}(\psi)$ has degree
$2n$ in $QH_*(M;\Lambda)$.

As pointed out in the introduction, at a first glimpse formula  (\ref{e:seidel}) of
Seidel's representation looks different from that of McDuff-Tolman \cite{mt}. When $(M,\omega)$
admits a Hamiltonian circle action $\psi$, then the definition of $\mathcal{S}(\psi)$ simplifies. 
For instance $P_\psi$ is isomorphic with the Borel quotient $S^3\times_{S^1}M$, where $S^3$ corresponds
to the total space of the Hopf fibration $S^3\to S^2$.  Hence if $p\in M$, is a fixed point
then the inclusion of $p$ into $M$, induces a section $\sigma_\textup{max}:S^2\to S^3\times_{S^1}M\simeq P_\psi$.
Thus there is a preferred section $\sigma_\textup{max}$ when  $(M,\omega)$
admits a Hamiltonian circle action. Further the index of the summation in Eq. (\ref{e:seidel}),
that is $A\in H^{sec}_2(P_{\psi})$ can be substituted by $\sigma_\textup{max}+B$ where
$B\in H_2(M,\mathbb{Z})$ is a spherical class. The rest of the details can be found
in Prop. 3.3 of \cite{mt}.
 
%%%%%%%%%%%%%%%%%%%%%%%%%%%%%%%%%%%%%%%%%%%%%%%%%%%%%%%%%%%%%%%
%%%%%%%%%%%%%%%%%%%%%%%%%%%%%%%%%%%%%%%%%%%%%%%%%%%%%%%%%%%%%%%
\section{Quantum product and the K\"unneth formula }
\label{s:quantkun}
%%%%%%%%%%%%%%%%%%%%%%%%%%%%%%%%%%%%%%%%%%%%%%%%%%%%%%%%%%%%%%%
%%%%%%%%%%%%%%%%%%%%%%%%%%%%%%%%%%%%%%%%%%%%%%%%%%%%%%%%%%%%%%%

So far we have only described the additive structure of quantum homology
$QH_*(M;\Lambda)=H_*(M)\otimes_\mathbb{Z}\Lambda$. However $QH_*(M;\Lambda)$ has the structure of a ring where
the operation is called quantum product. 

The quantum product is defined in terms of Gromov-Witten invariants, which are
a slide different from the ones discussed in the previous section. 
Consider $(M,\omega)$ a closed monotone symplectic manifold, homogeneous elements
$a,b,c\in H_*(M)$, and $A\in H_2(M)$ a spherical class. Then we have the
Gromov-Witten invariant
$
\textup{GW}^M_{A,3}(a,b,c),
$
which is the number of holomorphic curves that represent the class $A$
and intersect the cycles that represent the classes $a,b$ and $c$.
Here the degree of the homology classes must satisfy the equation
$
\textup{deg}(a)+\textup{deg}(b)+\textup{deg}(c)=4n-2c_1(A),
$
otherwise the invariant is defined as zero.
(See \cite{msjholo}, Ch. 7.) Now let $\{e_\nu\}_{\nu\in I}$ be a  base of the free $\mathbb{Z}$-module
$H_*(M)$, and $\{e^*_\nu\}_{\nu\in I}$ be the dual basis with respect to the intersection
product. That is $e^*_\nu\cdot e_\mu=\delta_{\nu,\mu}.$
Then if $a,b\in H_*(M)$ are homogeneous classes the quantum product 
$a*b$ is defined as
$$
a*b=\sum_{\nu\in I} \sum_A \textup{GW}^M_{A,3}(a,b,e_\nu)e^*_\nu
\otimes q^{-c_1(A)}t^{-\omega(A)},
$$
where $c_1$ is the first Chern class of $(M,\omega)$, and the sum
runs over all spherical classes $A\in H_2(M)$. Observe  that
$\textup{deg}(a*b)=\textup{deg}(a)+\textup{deg}(b)-2n$. Finally the quantum product
extends $\Lambda$-linearly to all $QH_*(M;\Lambda)$. Note that the identity element under
quantum multiplication corresponds
to the fundamental class $1=[M]\in H_{2n}(M)$.

An important fact about quantum homology is that 
the K\"unneth formula holds under a mild constraint. 
Let $(M,\omega)$ and $(N,\eta)$ be closed symplectic manifolds which are
 monotone with the \textit{same constant}. 
Thus the Gromov-Witten
invariants of $M\times N$ are well defined.
Let $a,b,c\in H_*(M\times N)$ be homogeneous classes such 
that the projections to $H_*(M)$ are denoted by $a_1,b_1$ and $c_1$, and  similarly
$a_2,b_2,c_2\in H_*(N)$. Then we have the following relation between
the Gromov-Witten invariants of $M, N$ and $M\times N$,
\begin{eqnarray}
\label{e:prodGW} 
\textup{GW}^{M\times N}_{A,3}(a,b,c)=
\textup{GW}^{M}_{A_1,3}(a_1,b_1,c_1)
\textup{GW}^{N}_{A_2,3}(a_2,b_2,c_2)
\end{eqnarray}
where $A\in H_2(M\times N)$ is a spherical class and $A_1,A_2$ correspond
to the projection of $A$ to $H_2(M)$ and $H_2(N)$ respectively. As a consequence 
we get the K\"unneth formula
for quantum homology
$$
QH_*(M\times N;\Lambda)\simeq
QH_*(M;\Lambda)\otimes_{\Lambda} QH_*(N;\Lambda).
$$
For more details see \cite{msjholo}.
With this at hand we  conclude that the maps $\kappa$ and $\kappa^\prime$
that appear in the main theorems are  ring homomorphisms under quantum multiplication.

\begin{lem} Let $(M,\omega)$ and $(N,\eta)$ be closed symplectic manifolds which 
are monotone with the same constant. Then the map of Thm. \ref{t:commdiag},
$$
\kappa:QH_*(M;\Lambda)\to QH_{*+\textup{dim}(N)}(M\times N;\Lambda)
$$
is a ring homomorphism under the quantum product.
\end{lem}
\begin{pf}
Let $\alpha_1,\alpha_2\in H_*(M)$. We must show that
$$
(\alpha_1\otimes [N])*(\alpha_2\otimes [N])
=(\alpha_1*\alpha_2)\otimes [N].
$$

This is a consequence of the K\"unneth formula for the quantum
homology ring. For
\begin{eqnarray*}
(\alpha_1\otimes [N])*(\alpha_2\otimes [N])
&=&(\alpha_1*\alpha_2)\otimes([N]*[N])\\
&=&(\alpha_1*\alpha_2)\otimes[N]
\end{eqnarray*}
since the fundamental class $[N]$ is the identity in $QH_*(N;\Lambda)$ under
quantum multiplication. Thus $\kappa$ is a ring homomorphism.
\end{pf}

A similar argument shows that the map $\kappa^\prime$
is a ring homomorphism.

\begin{lem} Let $(M,\omega)$ be a closed symplectic manifold such that
$\pi_2(M)$ is trivial. Then the map of Thm. \ref{t:2},
$$
\kappa^\prime:QH_*(M;\Lambda)\to QH_{*+*}(M\times M;\Lambda)
$$
is a ring homomorphism under the quantum product.
\end{lem}
\begin{pf}
By the K\"unneth formula, write the map $\kappa^\prime$
as
$$
\kappa^\prime:QH_*(M;\Lambda)\to QH_{*}(M;\Lambda)\otimes_\Lambda QH_{*}(M;\Lambda)
$$
where $\kappa^\prime(x)=x\otimes x$.  If follows again from
the K\"unneth formula  that for $x,y\in QH_{*}(M;\Lambda)$,
\begin{eqnarray*}
 \kappa^\prime(x*y)&=&(x*y)\otimes(x*y)\\
&=&(x\otimes x)*(y\otimes y)\\
&=&\kappa^\prime(x)*\kappa^\prime(y).
\end{eqnarray*}

Therefore $\kappa^\prime$ is a ring homomorphism under quantum multiplication.
\end{pf}

%%%%%%%%%%%%%%%%%%%%%%%%%%%%%%%%%%%%%%%%%%%%%%%%%%%%%%%%%%%%%%%%%%
%%%%%%%%%%%%%%%%%%%%%%%%%%%%%%%%%%%%%%%%%%%%%%%%%%%%%%%%%%%%%%%%%%
\section{Proof of the main result}
%%%%%%%%%%%%%%%%%%%%%%%%%%%%%%%%%%%%%%%%%%%%%%%%%%%%%%%%%%%%%%%%%%
%%%%%%%%%%%%%%%%%%%%%%%%%%%%%%%%%%%%%%%%%%%%%%%%%%%%%%%%%%%%%%%%%%

Let $(M,\omega)$ and $(N,\eta)$ be closed symplectic manifolds as in
Thm. \ref{t:commdiag}. That is $M$ is monotone and $\pi_2(N)=0$. Then
the product symplectic manifold $(M\times N,\omega\oplus\eta)$ is also
monotone, and therefore  Seidel's representation 
is well defined on  $(M\times N,\omega\oplus\eta)$.

\begin{lem}
\label{l:jacs}
Let $J\in \mathcal{J}(P_\psi,\pi,\Omega)$ where $\tilde\omega$
is a coupling form of $\pi:P_\psi\to S^2$ and $\Omega=
\tilde\omega + K\pi^*(\omega_0)$ as in Section \ref{s:smallqh}.
 If  $J^\prime$ is
a $\eta$-compatible almost complex structure on $TN$, then
$J\oplus J^\prime\in \mathcal{J}(P_\psi\times N,\pi_0, 
\tilde\omega\oplus\eta + K\pi_0^*(\omega_0))$.
\end{lem}
\begin{pf}
Let $p:P_\psi\times N\to P_\psi$ be the projection map. Then
$\pi\circ p=\pi_0$. On $P_\psi\times N$ we have the symplectic 
form $\Omega^\prime=\tilde\omega\oplus\eta +K\pi_0^*(\omega_0)$. 
Since $p^*\circ\pi^* =\pi_0^*$ then
$
\Omega^\prime=\Omega\oplus\eta
$
on $T(P_\psi\times N)\simeq TP_\psi\oplus TN$.
Thus if $J$ is an almost complex structure on $TP_\psi$ which
is $\Omega$-compatible and $J^\prime$  a $\eta$-compatible
 almost complex structure on $TN$, we have
\begin{eqnarray*}
\Omega^\prime(J\oplus J^\prime,J\oplus J^\prime)&=&
 \Omega(J,J)\oplus \eta( J^\prime,J^\prime)\\
&=&\Omega\oplus \eta\\
&=&\Omega^\prime.
\end{eqnarray*}
Hence $J\oplus J^\prime$ is an $\Omega^\prime$-compatible almost complex
structure on $T(P_\psi\times N).$

Assume that $J$ is such that the projection $\pi:(P_\psi,J)\to (S^2,j_0)$
is holomorphic. That is $d\pi\circ J=j_{0}\circ d\pi.$ Since
$d\pi_0=d\pi\circ dp$, we have
\begin{eqnarray*}
d\pi_0\circ (J\oplus J^\prime)&=&d\pi\circ dp\circ (J\oplus J^\prime)\\
&=&(d\pi\circ J) \oplus 0\\
&=& (j_{0}\circ d\pi)\oplus 0\\
&=& j_{0}\circ d\pi_0.
\end{eqnarray*}
Therefore, $J\oplus J^\prime\in \mathcal{J}(P_\psi\times N,\pi_0, 
\tilde\omega\oplus\eta + K\pi_0^*(\omega_0))$.
\end{pf}

The next proposition, is basically a restatement of Eq.
(\ref{e:prodGW}), but for the Gromov-Witten invariants that
are involved in the definition of the Seidel representation.

\begin{prop} 
\label{p:gweq}
Let $A\in H_2(P_{\psi\times \textup{id}_N}; \mathbb{Z})$ be a spherical class. 
Denote by $A_1:=(\lambda_1)_*(A)$ the induced
spherical class  in $H_2(P_\psi; \mathbb{Z})$. 
Then 
$$
\textup{GW}^{P_{\psi\times \textup{id}_N},\textup{\textbf{w}}}_{A,1}(\alpha\otimes
[\textup{pt}])=
\textup{GW}^{P_\psi,\textup{\textbf{w}}}_{A_1,1}(\alpha)
$$ 
for all $\alpha\in H_*(M)$.
\end{prop}
\begin{pf}
Let $\tilde\omega$ be a coupling form of $\pi:P_\psi\to S^2$ and  $J\in \mathcal{J}_{reg}(P_\psi,\pi,\Omega)$.
Thus
if $J^\prime$ is a $\eta$-compatible almost complex structure on $TN$, it follows from  
 Lemma \ref{l:jacs}, that
$J\oplus J^\prime\in \mathcal{J}(P_\psi\times N,\pi_0, \tilde\omega\oplus\eta +K\pi_0^*(\omega_0))$.
We must show that $J\oplus J^\prime$ is regular.

Let $u:S^2\to P_\psi\times N$ be a $(J\oplus J^\prime)$-holomorphic section that represents the
class $A\in H_2(P_\psi\times N; \mathbb{Z})$. Since $\pi_2(N)=0$, we may assume that
 $u=(u_0,q_0)$, where
$u_0:S^2\to P_\psi$ is a $J$-holomorphic section. Since $J$ is regular and 
$u_0$ is $J$-holomorphic, we know that
the linearized operator 
$$
D_{u_0}: \Omega^0(S^2,(u_0)^*(TP_\psi)) \to \Omega^{0,1}(S^2,(u_0)^*(TP_\psi))
$$
is onto. For the curve $u=(u_0,q_0)$, we have the linearized operator
$$
D_{u}: \Omega^0(S^2,(u_0)^*(TP_\psi)\oplus (S^2\times \mathbb{R}^{2m})) \to 
\Omega^{0,1}(S^2,(u_0)^*(TP_\psi)\oplus (S^2\times \mathbb{R}^{2m}))
$$
where $\textup{dim}(N)=2m.$ In this situation the operator
$D_u$ splits as the sum of $D_{u_0}$ and $\bar\partial$. See \cite{msjholo}, Rmk. 6.7.5. But the Cauchy-Riemann
operator $\bar\partial$ is also surjective,
 thus  it follows that $D_{u}$ is also onto. Therefore $J\oplus
J^\prime$ is regular.

Henceforth, $\mathcal{M}^\textbf{w}_1(A; J\oplus J^\prime)$ and $\mathcal{M}^\textbf{w}_1(A_1;J)$
are smooth oriented manifolds and can be use to compute the corresponding Gromov-Witten 
invariant. Let $\alpha\in H_*(M)$, since $\pi_2(N)=0$ we have the same intersection points
for the pseudocycle
$$
\iota^{-1}\circ ev_\textbf{w}:\mathcal{M}^\textbf{w}_1(A_1;J)\to M
$$
with $\alpha$; and the pseudocycle
$$
\iota_0^{-1}\circ ev_\textbf{w}:\mathcal{M}^\textbf{w}_1(A; J\oplus J^\prime)\to M\times N.
$$
with $\alpha\otimes[\textup{pt}]$. Hence
$
\textup{GW}^{P_{\psi\times \textup{id}_N},\textup{\textbf{w}}}_{A,1}(\alpha\otimes
[\textup{pt}])=
\textup{GW}^{P_\psi,\textup{\textbf{w}}}_{A_1,1}(\alpha).
$
\end{pf}

Now by Prop. \ref{p:gweq}, there is a similar relation between  the 
homology classes $\mathcal{S}_A(\psi)$
and  $\mathcal{S}_A(\psi\times \textup{id}_N)$. 

\begin{lem}
\label{l:sames}

Let $A\in H_2(P_{\psi\times \textup{id}_N}; \mathbb{Z})$ and 
$A_1\in H_2(P_\psi; \mathbb{Z})$ as in Prop. \ref{p:gweq}. Then the identity
$$
\mathcal{S}_A(\psi\times\textup{id}_N)=\mathcal{S}_{A_1}(\psi)\otimes[N]
$$
holds in $H_*(M\times N)$.
\end{lem}

\begin{pf}
Let $\alpha\in H_*(M)$ and $\beta\in H_*(N)$  such that the sum of the degrees of 
$\alpha$ and $\beta$ is $-2c_\psi(A_1)$. Then by the definition of the invariant
$\textup{GW}^{P_\psi,\textbf{w}}_{A_1,1}$,
\begin{eqnarray}
\label{eqgw1}
(\mathcal{S}_{A_1}(\psi)\otimes[N]) \cdot_{M\times N}(\alpha\otimes \beta)&=&
(\mathcal{S}_{A_1}(\psi)\cdot_M \alpha)\cdot([N] \cdot_{N}\beta)\nonumber\\
 &=&\textup{GW}^{P_\psi,\textbf{w}}_{A_1,1}(\alpha) \cdot([N] \cdot_{N}\beta).
\end{eqnarray}
The terms in this equation
 are all equal to zero unless $\alpha$ has degree $-2c_\psi(A_1)$
and $\beta$ has degree $0$, that is $\beta=[\textup{pt}]$. 
In this case, note that $[N]\cdot_{N}[\textup{pt}]=1$. 

On the other hand by the definition of $\textup{GW}^{P_{\psi\times\textup{id}_N},\textbf{w}}_{A,1}$
we have
\begin{eqnarray}
\label{eqgw2}
\mathcal{S}_A(\psi\times \textup{id}_N) \cdot_{M\times N}(\alpha\otimes \beta) & =&
\textup{GW}^{P_{\psi\times\textup{id}_N},\textbf{w}}_{A,1}(\alpha\otimes \beta).
\end{eqnarray}
Since the class $\alpha$ is in $H_*(M)$, then by definition of the Gromov-Witten
invariant, it follows that
$\textup{GW}^{P_{\psi\times\textup{id}_N},\textbf{w}}_{A,1}(\alpha\otimes \beta)$ is
zero unless $\beta=[\textup{pt}]$. Hence from Prop. \ref{p:gweq},  Eqs.
(\ref{eqgw1}) and (\ref{eqgw2}) are equal. That is,
\begin{eqnarray*}
(\mathcal{S}_{A_1}(\psi)\otimes[N]) \cdot_{M\times N}(\alpha\otimes \beta)
& =&
\mathcal{S}_A(\psi\times \textup{id}_N) \cdot_{M\times N}(\alpha\otimes \beta) 
\end{eqnarray*}
for all $\alpha\in H_*(M)$ and $\beta\in H_*(N)$. Therefore 
$\mathcal{S}_{A_1}(\psi)\otimes[N]=\mathcal{S}_A(\psi\times \textup{id}_N)$.
\end{pf}

\noindent \textbf{Proof of Thm. \ref{t:commdiag}.}
First of all, since $\pi_2(N)$ is trivial we have that 
$$
(\lambda_1)_*: H_2(P_{\psi\times \textup{id}_N};\mathbb{Z})
\to H_2(P_{\psi};\mathbb{Z})
$$
induces a one-to-one correspondence between the
section classes of $P_\psi$ and the section classes of $P_{\psi\times\textup{id}_N}$.
That is, $H^{sec}_2(P_{\psi})\simeq
H^{sec}_2(P_{\psi\times \textup{id}_N})$.
Hence the sum on the definition   of the elements
$\mathcal{S}(\psi)$ and $\mathcal{S}(\psi\times\textup{id}_N)$ is
defined over the same set. 

For $A\in H^{sec}_2(P_{\psi\times \textup{id}_N};\mathbb{Z})$, we have that
$(\lambda_2)_*(A)=0$ since $\pi_2(N)=0$. Therefore
$$
c_{\psi\times \textup{id}_N}(A)=c_\psi(A_1)
\mbox{\hskip .5cm and \hskip .5cm }
\tilde\omega_{\psi\times \textup{id}_N}(A)=
\tilde\omega_\psi(A_1)
$$
by  Eqs. (\ref{eq.samec}) and (\ref{eq.samew}).
Finally from  Lemma \ref{l:sames} we have
$\mathcal{S}_A(\psi\times\textup{id}_N)=\mathcal{S}_{A_1}(\psi)\otimes[N]$.
Therefore $\mathcal{S}(\psi\times\textup{id}_N)=
\mathcal{S}(\psi)\otimes[N]$. 
\hfill$\Box$ \medskip

\medskip
Consider the group homomorphism 
$
\tau_0:\pi_1(\textup{Ham}(M,\omega))
\to \pi_1(\textup{Ham}(M\times M,\omega\oplus\omega))
$
defined as $\tau_0(\psi)=\textup{id}_M\times \psi.$ Define
the map
$$
\kappa_0:QH_*(M;\Lambda)\to QH_{*+\textup{dim}(M)}(M\times M;\Lambda)
$$
on homogeneous elements by
$\kappa_0(\alpha\otimes q^rt^s)=(\alpha\otimes [M]) q^rt^s$, where $\alpha\in H_*(M)$,
and extend it $\Lambda$-linearly to all $QH_*(M;\Lambda)$.
Then as in Thm \ref{t:commdiag}, we have that
\begin{eqnarray}
\label{e:comnmu}
\mathcal{S}\circ\tau_0=\kappa_0\circ\mathcal{S}.
\end{eqnarray}

Then Thm. \ref{t:commdiag} together with Eq. (\ref{e:comnmu}) provide a
a proof of Thm. \ref{t:2},

\noindent \textbf{Proof of Thm. \ref{t:2}.} 
Observe that $\tau^\prime(\psi)=\tau(\psi)\circ \tau_0(\psi)$. Then
applying Seidel's representation we get  from Thm. \ref{t:commdiag}
and Eq. (\ref{e:comnmu}) that 
\begin{eqnarray*}
\mathcal{S}\circ \tau^\prime(\psi)
&= & \mathcal{S}(\tau(\psi)\circ \tau_0(\psi)) \\
&= & \mathcal{S}(\tau(\psi))*\mathcal{S}( \tau_0(\psi)) \\
&= & (\kappa\circ\mathcal{S}(\psi))*(\kappa_0\circ\mathcal{S}(\psi))\\
&= & (\mathcal{S}(\psi)\otimes [M])*([M]\otimes\mathcal{S}(\psi)).
\end{eqnarray*}

Then by the K\"unneth formula and the fact that
$[M]$ is the indentity on $QH_*(M;\Lambda)$ we get
\begin{eqnarray*}
\mathcal{S}\circ \tau^\prime(\psi)
&= & (\mathcal{S}(\psi)\otimes [M])*([M]\otimes\mathcal{S}(\psi))\\
&= & (\mathcal{S}(\psi)*[M])\otimes([M]*\mathcal{S}(\psi))\\
&= & \mathcal{S}(\psi)\otimes\mathcal{S}(\psi)\\
&= & \kappa^\prime(\mathcal{S}(\psi)).
\end{eqnarray*}
\hfill$\Box$ \medskip

%\ack

\end{document}